\documentclass{article}
\usepackage{amssymb}
\usepackage{latexsym}
\usepackage{amsfonts}
\usepackage{amsbsy}
\usepackage{subfigure}
\usepackage{mathtools}
\usepackage{graphicx}
\usepackage{algorithm}
\usepackage{algpseudocode}
\usepackage{caption}
\usepackage{bbm}
\usepackage{booktabs}
\usepackage{cite}

\usepackage{amsthm}

\usepackage{url}
\usepackage{xcolor}
\definecolor{newcolor}{rgb}{.8,.349,.1}

\graphicspath{{figures/}}

\theoremstyle{plain} 
\newtheorem{remark}{Remark}
\newcommand{\opI}{\mathcal{I}}
\newcommand{\opD}{\mathcal{D}} 

\newcommand{\opS}{\mathcal{S}} 
 
\newcommand{\Ord}{\mathcal{O}}

\let\emptyset\varnothing

\newcommand{\N}{\mathbb{N}}
\newcommand{\R}{\mathbb{R}}
\newcommand{\C}{\mathbb{C}}
\newcommand{\Sfar}{\mathbb{S}^2_N}
\newcommand{\Pnear}{\mathbb{P}_N}
\newcommand{\Real}[1]{\text{Re}\left\{#1\right\}}

\newcommand{\ome}{\omega}

\newcommand{\ii}{\text{i}}
\newcommand{\sgn}{\text{sgn}} 

\newcommand{\farErr}{\varepsilon_\textit{far}}
\newcommand{\nearErr}{\varepsilon_\textit{near}}

\newcommand{\del}{\partial}
\newcommand{\Dom}{\Omega}
\newcommand{\Bdy}{\Gamma}

\newcommand{\dist}[1]{\text{dist}\left(#1\right)}

\newcommand{\norm}[1]{\left\Vert#1\right\Vert}

\newcommand{\floor}[1]{\left\lfloor#1\right\rfloor}

\newcommand{\bgamma}{\boldsymbol\gamma}
\newcommand{\ykd}{y^{\mathbf{k}}_d}
\newcommand{\Bkd}{B^d_{\mathbf{k}}}
\newcommand{\Wkd}{W^{\mathbf{k},d}_j(s,\theta,\phi,y)}

\DeclareMathOperator*{\argmin}{arg\,min}

\title{IFGF-accelerated high-order integral equation solver\\ for
  acoustic wave scattering}

\author{Edwin Jimenez$^*$, Christoph Bauinger\footnote{Computing and Mathematical Sciences, Caltech, Pasadena, CA 91125, USA} ~and Oscar P. Bruno$^*$\footnote{Corresponding author.}}

\begin{document}

\maketitle

\begin{abstract}
  We present an accelerated iterative boundary-integral solver for the
  numerical solution of problems of time-harmonic acoustic scattering
  by general surfaces in three-dimensional space. The proposed method
  relies on the recently introduced high-order rectangular-polar
  algorithm (RP) for evaluation of singular integrals, and it utilizes
  an extended version of the novel IFGF (Interpolated Factored Green
  Function) acceleration scheme---demonstrating, in particular, the
  first application of the IFGF approach as part of a full scattering
  solver.  Exploiting slow variations in certain factored forms of the
  scattering Green function, the IFGF algorithm enables the fast
  evaluation of fields generated by groups of sources on the basis of
  a recursive interpolation scheme.  Relying on Chebyshev expansions
  and regularizing changes of variables, in turn, the RP algorithm
  accurately evaluates challenging singular and near-singular
  Green-function interactions.  A parallel OpenMP implementation of
  the overall algorithm is presented, and numerical experiments
  confirm that the expected accuracy, and the overall $\Ord(N\log N)$
  computational cost as the frequency and discretization sizes are
  increased, are observed in practice.  Numerical examples include
  acoustic scattering by spheres of up to $128$ wavelengths in
  diameter, an $80$-wavelength submarine, and a turbofan nacelle that
  is more than $80$ wavelengths in size, requiring, on a 28-core
  computer, computing times of the order of a few minutes per
  iteration and a few tens of iterations of the GMRES iterative
  solver.\looseness=-1
\end{abstract}

\noindent {\em Keywords: Integral equations, Fast solver, Scattering
  problem, High-order accuracy, General geometries}



\section{Introduction}

This paper presents an accelerated iterative boundary-integral solver
for the numerical solution of problems of time-harmonic acoustic
scattering by general surfaces in three-dimensional space. On the
basis of the iterative linear-algebra solver GMRES~\cite{Saad1986},
the proposed algorithm relies on the high-order rectangular-polar (RP)
singular integration method~\cite{BrunoGarza2020}, and it utilizes an
extended version of the recently-introduced IFGF acceleration
scheme~\cite{BauingerBruno2021} that is applicable to the
combined-field acoustic integral
equations~\cite{ColtonKress2013IEMethods}---providing, in fact, the
first application of the IFGF acceleration approach as part of a full
scattering solver. An OpenMP parallel implementation of the algorithm
is demonstrated, which can produce accurate solutions for complex,
acoustically-large engineering structures in computing times of a few
minutes per iteration in a relatively small ($28$ core) computer
system.

As is well known, boundary integral methods provide a number of
advantages (notably, they only require discretization of the
scattering surfaces, they do not suffer from the dispersion or
pollution errors associated with differential formalisms, and they
inherently enforce the condition of radiation at infinity). But they
do give rise to certain challenges---concerning the accurate
evaluation of the associated singular integrals, and the $\Ord(N^2)$
computational operator-evaluation cost, for an $N$-point surface grid,
that results from straightforward integral-equation computational
implementations.  Like the Fast Multipole Method (FMM) and other
acceleration approaches~\cite{Bleszynski1996,BrunoKunyansky2001JCP,
  ChengEtAl2006,EngquistYing2007,
  MichielssenBoag1996,PhillipsWhite1997,PoulsonEtAl2014,
  Rokhlin1993,YingEtAl2003} that rely on the Fast Fourier Transform
(FFT), the aforementioned IFGF strategy reduces the computing cost to
$\Ord(N\log N)$ operations. The IFGF algorithm, however, is not based
on previously-employed acceleration elements such as the Fast Fourier
Transform (FFT), special-function expansions, high-dimensional
linear-algebra factorizations, translation operators, equivalent
sources, or parabolic scaling.  Instead, the IFGF method relies on an
interpolation scheme of factored forms of the operator kernels which,
when applied recursively to larger and larger groups of Green function
sources, gives rise to the desired $\Ord(N\log N)$ accelerated
evaluation---in a manner that lends itself to effective
parallelization under both shared-memory systems and large
distributed-memory hardware
infrastructures~\cite{BauingerBruno_parallel_2021}.

The narrative description of the IFGF algorithm presented in
Section~\ref{sec:IFGF} below, which extends and supplements the
description in the previous contribution~\cite{BauingerBruno2021},
includes special considerations concerning the application of the IFGF
method to the aforementioned combined field equation; it presents two
different strategies for IFGF acceleration of the corresponding
double-layer operator, each one advantageous under complementary
structural settings; and it emphasizes an intuitive presentation of
the necessary algorithmic details.  The RP algorithm, together with
the proposed strategy for the integration of the RP and IFGF
algorithms, in turn, are presented in Section~\ref{sec:RecPolar}.

The RP approach utilizes a CAD-like (Computer Aided Design)
representation of the scattering surface as the union of structured
non-overlapping surface patches.  Chebyshev expansions of the
underlying integral densities are used on each patch, and the
necessary singular integral operators are produced with high-order
accuracy by means of suitable singularity-smoothing changes of
variables. The RP method, which yields high order accuracy in spite of
the Green function singularity, can, without IFGF acceleration, be
utilized as the sole basis of a complete high-order scattering
solver. But its $\Ord(N^2)$ computational cost is prohibitive for
acoustically-large problems.

In the combined IFGF/RP approach, the relatively few (but challenging)
singular and near-singular Green function interactions are obtained on
the basis of the RP method, and the numerous non-singular
contributions are incorporated, in an accelerated manner, by means of
the IFGF algorithm. In view of the $\Ord(N)$ computational cost
required for evaluation of singular and near-singular interactions by
the RP algorithm, an overall $\Ord(N\log N)$ complexity for the
combined IFGF/RP algorithm results. The description of the integrated
IFGF/RP algorithm presented in this paper may additionally be used as
a blueprint for use of IFGF acceleration in conjunction with other
integral equation discretization strategies, such as the Method of
Moments~\cite{harrington1992} and other Galerkin and Nystr\"om
boundary integral discretization approaches. The combined IFGF/RP
methodology is demonstrated in this paper by means of an OpenMP
parallel implementation suitable for shared-memory computer systems;
the development of related acoustic and electromagnetic scattering
solvers on large distributed-memory systems is left for future work.

A variety of numerical examples presented in Section~\ref{sec:NumRes}
show that, as suggested above, the proposed approach leads to the
efficient solution of large scattering problems for complex geometries
on small parallel computers.  Numerical examples include acoustic
scattering by spheres up to $128$ wavelengths in diameter, an
$80$-wavelength submarine, and a turbofan nacelle that is more than
$80$ wavelengths in size, each one requiring, on a 28-core computer,
computing times of the order of a few minutes per iteration and a few
tens of iterations of the GMRES iterative solver. Convergence studies
were used to determine the accuracy of the solution for the
engineering structures considered, for which exact solutions are, of
course, not available. In particular, the computational results
obtained confirm the validity in practice of the theoretical
$\Ord(N\log N)$ computing-cost estimate as the frequency and
discretization are increased. Comparisons of the proposed solver with
the recent FMM-based
algorithms~\cite{WalaKlockner2019,AbduljabbarEtAl2019} are provided in
Section~\ref{compari}.
 
This paper is organized as follows. Preliminary topics, including the
integral-equation formulations utilized and the proposed
non-overlapping patch surface representation used, are briefly
reviewed in Section~\ref{sec:IntEqn}.  Section~\ref{sec:IFGF} then
presents the IFGF algorithm, with a focus on the extensions needed for
IFGF-acceleration of the double-layer integral operator, and including
details concerning IFGF OpenMP-based parallelization.  The proposed
high-order RP algorithm for evaluation of singular integral operators,
and the integration of the RP and IFGF algorithms, are then described
in Section~\ref{sec:RecPolar}. The numerical results presented in
Section~\ref{sec:NumRes} demonstrate the character of the overall
IFGF-based solver by means of a variety of numerical experiments.
Section~\ref{sec:Conclusion}, finally, presents a few concluding
remarks.

\section{Preliminaries}\label{sec:IntEqn}

\subsection{Scattering boundary-value problem}\label{sec:AcousticBVP}
We consider wave propagation in a homogeneous isotropic medium with
density $\rho$, speed of sound $c$, and in absence of
damping~\cite{ColtonKress2013IEMethods}.  Scattering obstacles are
represented by the closure of a bounded set $\Dom \subset \R^3$; the
complement $\bar{\Dom}^c$ of the closure $\bar{\Dom} \subset \R^3$ is
the propagation domain. For time-harmonic acoustic waves, the wave
motion can be obtained from the velocity potential
$U(x,t) = \Real{u(x) e^{-i\ome t}}$, where $\ome > 0$ is the angular
frequency, and where the complex-valued function $u(x)$ satisfies the
Helmholtz equation
\begin{equation}\label{eq:ScalarHelmholtz} 
  \Delta u(x) + k^2 u(x) = 0, \qquad x \in \R^3 \setminus \bar{\Dom}, 
\end{equation} 
with $k = \ome / c$; the corresponding acoustic wavelength is given by
$\lambda = 2\pi/k$. Denoting the boundary of $\Dom$ by $\Bdy$, the
\emph{sound-soft} obstacle case that we consider requires that $u = 0$
on $\Bdy$. Writing the total field $u(x) = u^i(x) + u^s(x)$, where
$u^i(x)$ is a given incident field which also satisfies Helmholtz
equation in a neighborhood of $\bar{\Dom}$, leads to an exterior
Dirichlet boundary-value problem for the scattered field $u^s(x)$:
\begin{align}\label{eq:AcousticBVP}
  \begin{cases}
  \Delta u^s(x) + k^2 u^s(x) = 0, &x \in \R^3 \setminus \bar{\Dom}, \\
  u^s(x) = -u^i(x),               &x \in \Bdy, \\
  |x| \left( \displaystyle\frac{x}{|x|} \cdot \nabla u^s(x) - i k u^s(x) \right) = 0, & |x| \to \infty.
  \end{cases}
\end{align}

\subsection{Integral representations and integral equations}\label{sec:IntRep}
To tackle the scattering problem~\eqref{eq:AcousticBVP} we express the
scattered $u^s$ as a combined-field
potential~\cite{ColtonKress2013IEMethods}
\begin{equation}\label{eq:CombinedLayer}
  u^s(x) = \int_{\Bdy} \left\{ \frac{\del \Phi(x,y)}{\del \nu(y)}  
                             - \ii \gamma \Phi(x,y) \right\} \varphi(y) \, dS(y), 
    \qquad x \in \R^3 \setminus \bar{\Dom}
\end{equation}
for a real \emph{coupling parameter} $\gamma \neq 0$, where $\nu$
denotes the outer normal to $\Gamma$, and where
\begin{equation}\label{eq:GreenFunction}
  \Phi( x, y ) = \frac{1}{4\pi} \frac{e^{\ii k |x-y|}}{|x-y|}
\end{equation}
denotes the corresponding Green function. (Equation~\eqref{eq:gamma} in the
numerical results section and associated discussion concern the impact of the
value of the parameter $\gamma$ on the performance of the numerical method.)
Noting that
\begin{equation}\label{double_layer}
 \frac{\partial \Phi(x,y)}{\partial \nu(y)} = \frac{e^{\ii k |x-y|}}{4\pi |x-y|}\left( 1 - \ii k |x-y| \right) 
 \frac{\langle x - y, \nu(y) \rangle}{|x-y|^2}
 = e^{\ii k |x-y|}\frac{\langle x - y, \nu(y) \rangle}{4\pi |x-y|^3}  - \ii ke^{\ii k |x-y|}\frac{\langle x - y, \nu(y) \rangle}{4\pi |x-y|^2},
\end{equation}
and letting
\begin{equation}\label{eq:scalarop}
  \opS [ \varphi ](x) = \int_{\Bdy} \Phi(x,y) \varphi(y) \, dS(y) 
  \quad   \mbox{and}\quad
  \opD [ \varphi ](x) = \int_{\Bdy} 
  \frac{\del \Phi(x,y)}{\del \nu(y)} \varphi(y) \, dS(y), \quad (x\in\Gamma),
\end{equation}
denote the single- and double-layer integral operators, respectively,
the unknown scalar density function $\varphi$ is obtained as the
solution of the combined-field integral equation
\begin{equation}\label{eq:CombLayIE}
  \frac{1}{2} \varphi(x) + \opD [\varphi](x) - \ii \gamma \opS[\varphi](x) = -u^i(x),
    \qquad x \in \Bdy,
\end{equation}
which guarantees that the boundary values on $\Gamma$, as required
by~\eqref{eq:AcousticBVP}, are attained. Once the density $\varphi$ has
been obtained, the scattered field $u^s$ everywhere outside $\Gamma$
can be obtained by direct integration from~\eqref{eq:CombinedLayer}.

For reference we define the functions
\begin{equation}\label{single_double_terms}
  V_1(x,y)= \frac{e^{\ii k |x-y|}}{|x-y|}, \quad V_2(x,y)= e^{\ii k |x-y|}\frac{\langle x - y, \nu(y) \rangle}{4\pi |x-y|^3},\quad\mbox{and}  \quad V_3(x,y)= e^{\ii k |x-y|}\frac{\langle x - y, \nu(y) \rangle}{4\pi |x-y|^2};
\end{equation}
in view of~\eqref{eq:GreenFunction} and~\eqref{double_layer}, we have
\begin{equation}\label{single_double_terms_2}
 \Phi( x, y )= \frac{1}{4\pi} V_1(x,y)\quad \mbox{and}\quad \frac{\partial \Phi(x,y)}{\partial \nu(y)} = \frac{1}{4\pi}V_2(x,y) - \frac{\ii k}{4\pi}  V_3(x,y).
\end{equation}

\subsection{Surface representation}\label{sec:SurfRep}
\begin{figure}[H]
    \centering
    \includegraphics[width=0.9
   \textwidth]{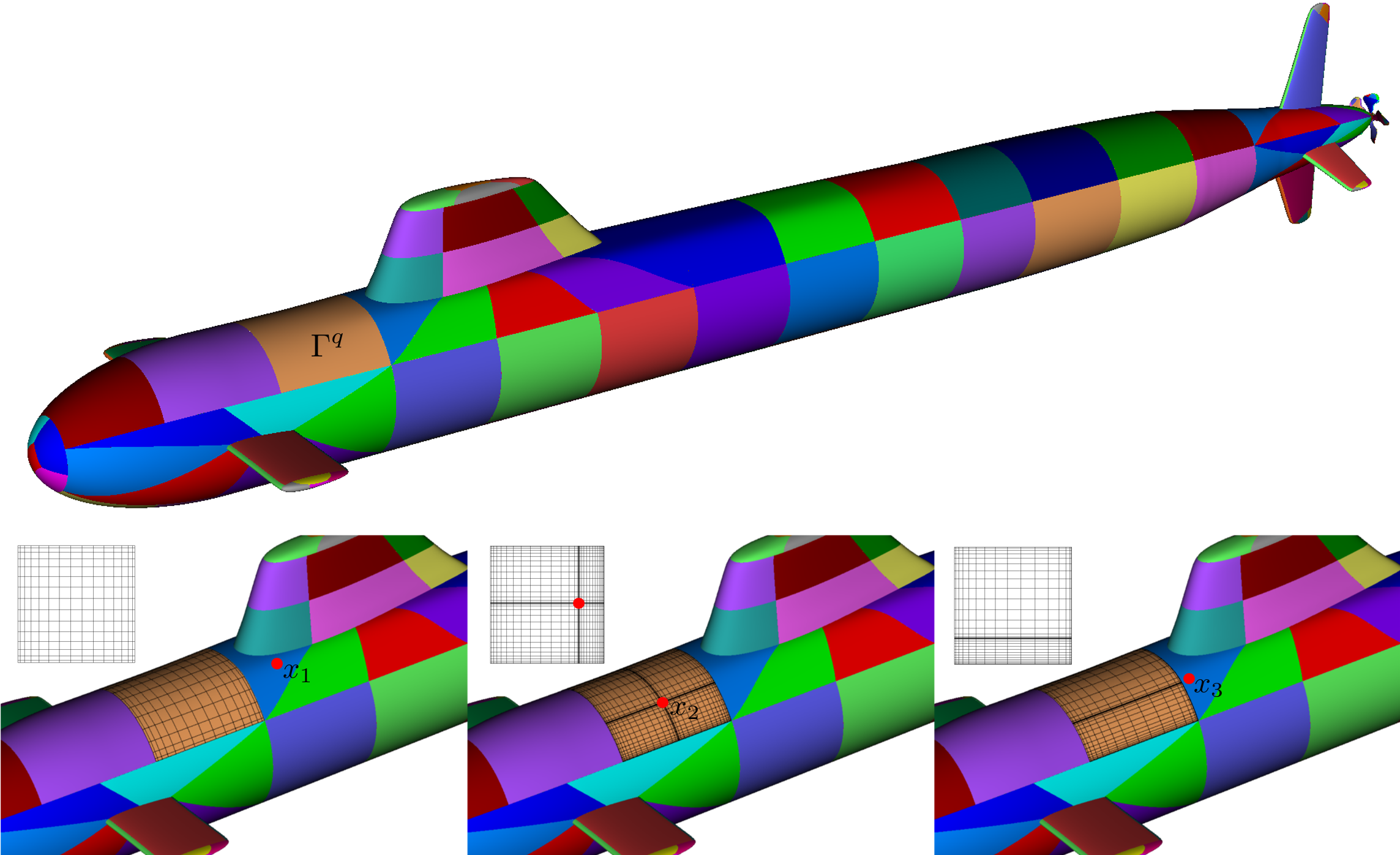}
   \caption{\small Upper panel: Non-overlapping multi-patch surface
     representation for a submarine model, with a patch denoted by
     $\Gamma^q$ depicted in light-brown color. Lower left panel: a
     patch $\Gamma^q$ overlaid with a Chebyshev discretization mesh
     induced by the parameter grid displayed in the leftmost inset
     parameter patch $(-1,1)^2 = (-1,1)\times (-1,1)$. An observation
     point $x_1$ far from the patch is shown in red, for which the
     Fej\'er first quadrature rule on the aforementioned parameter
     grid gives rise to accurate integration.  Lower center panel:
     same as the lower-left panel but for an evaluation point $x_2$
     within the patch $\Gamma^q$, and associated graded meshes in the
     parameter and physical patches. Lower right panel: same as the
     lower-center panel, for an evaluation point outside but near the
     patch $\Gamma^q$, and associated graded meshes in the parameter
     and physical patches.  }
    \label{fig:PatchesAndGrids}
\end{figure}

In order to obtain a computational implementation of
equation~\eqref{eq:CombLayIE}, we partition the scattering surface
$\Gamma$ as the disjoint union of a set of non-overlapping
parametrized component ``patches''---as illustrated in
Figure~\ref{fig:PatchesAndGrids} for the particular case of a
submarine geometry.  The surface patches we use are ``logical
quadrilaterals'' (LQ): each one of them is the image of a rectangular
reference domain under an certain vector-valued smooth parametrization
function---in such a way that physical corners and edges of the
surface $\Gamma$ only arise at points and curves at which smooth
patches meet.

In detail, for a given scattering surface $\Bdy$ we utilize a certain number
$Q$ of smooth parametrizations
\begin{equation}\label{param_def}
  y^q : (-1,1)^2 \to \R^3, \qquad (q=1,\dotsc,Q),
\end{equation}
each one of which maps a $uv$ reference domain $(-1,1)^2= (-1,1)\times (-1,1)$
onto an LQ patch
\begin{equation}\label{eq:GammaCover}
  \Gamma^q = y^q((-1,1)^2)\subset \R^3, \quad \text{so that, in all} \quad
  \Gamma = \bigcup_{q=1}^Q \Gamma^q.
\end{equation}
On the basis of such a geometric representation of the surface
$\Gamma$, a general surface integral operator
\begin{equation}
  \label{eq:gen-int-op}
  \opI [\varphi](x) = \int_{\Bdy} K(x,y) \varphi(y) \, dS(y),\quad x\in\Gamma,
\end{equation}
such as those in equation~\eqref{eq:scalarop}, can be evaluated in
patch-wise fashion,
\begin{equation}
  \label{eq:patch-int}
  \opI [\varphi](x) = \sum_{q=1}^Q \opI^q[\varphi](x)\quad \mbox{where}\quad \opI^q [\varphi](x) = \int_{\Bdy_q} K(x,y) \varphi(y) \, dS(y),\quad x\in\Gamma,
\end{equation}
and, thus, the evaluation of $\opI [\varphi](x)$ for $x\in\Gamma$ reduces to
evaluation of the patch operators $\opI^q [\varphi](x)$ for all $x\in\Gamma$
and $q = 1,\dotsc, Q$. The problem of evaluation of the patch operators $\opI^q
[\varphi](x)$, whose treatment require careful consideration of the position of
the point $x\in\Gamma$ relative to the patch $\Gamma_q$, is taken up in
Section~\ref{sec:RecPolar}.

\section{IFGF acceleration of discrete integral operators}\label{sec:IFGF}

As detailed in Section~\ref{sec:RecPolar}, the iterative numerical
solution of equation~\eqref{eq:CombLayIE} requires evaluation of
discrete versions of the integral operators $\opD$ and $\opS$
introduced in the previous section, for given numerical approximations
of the unknown density $\varphi$, on a suitably selected surface mesh
$\Bdy_N = \{x_1,\dots,x_N\}$, where $x_j$ ($1\leq j\leq N$) are
pairwise different surface discretization points. Per the discussion
in Section~\ref{sec:RecPolar}, the evaluation of such discrete
operators entails two main challenges: (a)~The accurate evaluation of
singular and near-singular contributions (namely, contributions to the
integrals in~\eqref{eq:scalarop} arising from values of $y$ close to a
given $x$, for which the singularity of the Green function must be
suitably considered); and (b)~The efficient evaluation of the vast
number of non-singular contributions---which, as discussed in
Section~\ref{sec:RecPolar}, reduces to evaluation of discrete
operators of the form
\begin{equation} \label{eq:field1} I_K(x_\ell) = \sum
  \limits_{\substack{m = 1 \\ m \neq \ell}}^N a_m K(x_\ell, x_m)
  ,\quad \ell = 1, \ldots, N,
\end{equation}
where $a_m \in \C$ denote complex numbers, and where $K$ denotes
either the Green function itself ($K(x,y)=\Phi( x, y )$) or its normal
derivative ($K(x,y)=\frac{\partial \Phi}{\partial \nu(y)}$).  The
present section addresses the challenges arising from point~(b): on
the basis of the IFGF recursive interpolation strategy described in
Sections~\ref{sec:GenIFGF} through~\ref{sec:CombLayIFGF} below, the
proposed approach effectively reduces the $\Ord{(N^2)}$ cost required
by direct evaluation of discrete operators of the
form~\eqref{eq:field1} to $\Ord{(N\log N)}$ operations,
in a manner
that remains efficiently parallelizable for challenging high-frequency
configurations.  The IFGF acceleration strategy was introduced for the
single layer kernel $K(x,y)=\Phi( x, y )$ in~\cite{BauingerBruno2021}
and is extended in in what follows to the double-layer kernel
$K(x,y)=\frac{\partial \Phi}{\partial \nu(y)}$.

\subsection{Factored Green function: analytic properties}\label{sec:GenIFGF}

\begin{figure}[h!]
  \centering
  \includegraphics[width=\textwidth]{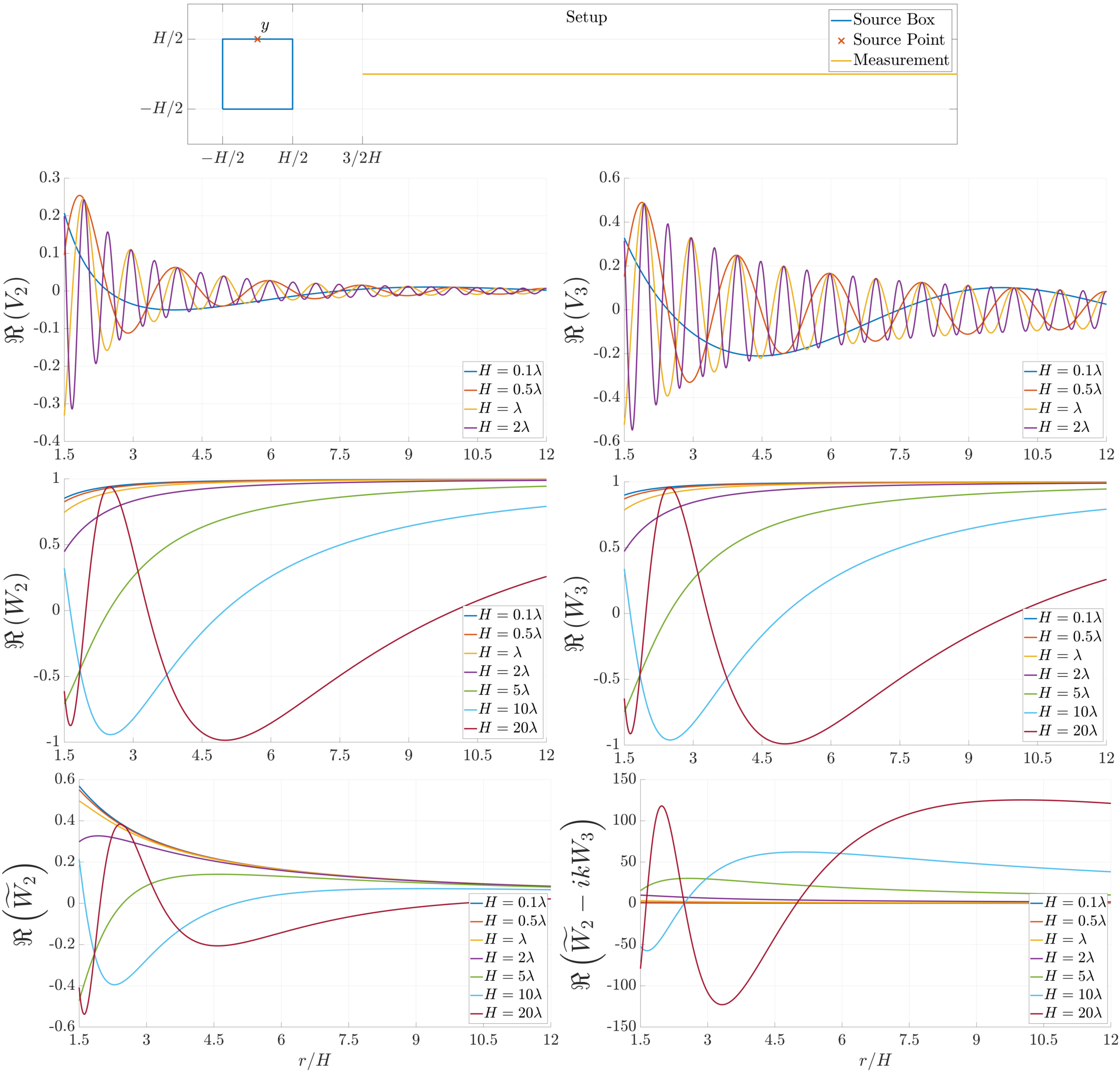}
  \caption{\small Source-point factorization test, set up as
    illustrated in the ``Setup'' panel. The subsequent two panels
    demonstrate the highly oscillatory character of the double-layer
    kernel components $V_2$ and $V_3$ (the graphs display the real
    parts $\Re (V_2)$ and $\Re (V_3)$ of $V_2$ and $V_3$), even for
    small values of $H$ ($H\leq 2\lambda$). The final four panels
    illustrate the corresponding slowly oscillatory character of
    various IFGF analytic factors introduced in
    Section~\ref{sec:GenIFGF}, for $H$ values up to $H=20\lambda$. The
    factor $W_1$, which is considered in~\cite{BauingerBruno2021}, and
    which is quite similar in form and character to $W_3$, is not
    included in this figure, for conciseness.}
  \label{fig:FactUnfaDL}
\end{figure}

We seek to efficiently evaluate the sum $I_K$ in~\eqref{eq:field1}
with kernel function $K$ equal to either the single- or double-layer
kernel (equations~\eqref{eq:GreenFunction} and~\eqref{double_layer},
respectively; cf.~\eqref{single_double_terms}
and~\eqref{single_double_terms_2}). In the case of the single-layer
kernel $\Phi( x, y )= \frac{1}{4\pi} V_1(x,y)$ we utilize the factored
expression
\begin{equation}\label{hgns}
   V_1(x,y) =\frac{e^{\ii kr}}{r}\frac{e^{\ii kr\left( \left |\frac{x- y^0}{r}-\frac{y- y^0}{r}\right| - 1\right)}}{\left |\frac{x- y^0}{r}-\frac{y- y^0}{r}\right|},\quad r = |x - y^0|,
\end{equation}
which, for points $y$ within a cubic box $B(y^0,H)$ of side $H$ centered at the
point $y^0 = (y^0_1,y^0_2,y^0_3)$,
\begin{equation} \label{eq:box} 
  B(y^0,H) \coloneqq \bigg[ y^0_1 - \frac{H}{2}, y^0_1 + \frac{H}{2} \bigg)
              \times \bigg[ y^0_2 - \frac{H}{2}, y^0_2 + \frac{H}{2} \bigg)
              \times \bigg[ y^0_3 - \frac{H}{2}, y^0_3 + \frac{H}{2} \bigg),
\end{equation}
and using a parameter $h$ and a new variable $s$ given by
\begin{equation} \label{eq:def_eps} h = \max \limits_{y \in B_{y^0}} |y- y^0| =
  \frac{\sqrt{3}}{2}H\quad\mbox{and}\quad s=\frac{h}{r},
\end{equation}
we re-express in the form
\begin{equation}\label{hgn2}
  V_1(x,y) =\frac{e^{\ii kr}}{r}
  \frac{e^{\ii k\frac{h}{s}\left(\left |\frac{x- y^0}{r}-\frac{y- y^0}{h}s\right| - 1\right)}}{ \left |\frac{x- y^0}{r}-\frac{y- y^0}{h}s\right|} = \frac{e^{\ii kr}}{r}\widehat{W}_1(x,y;y^0).
\end{equation}
As indicated in what follows, for each fixed $y$, the second factor
\begin{equation}\label{V3}
  \widehat{W}_1(x,y;y^0)= \frac{e^{\ii k\frac{h}{s}\left(\left |\frac{x- y^0}{r}-\frac{y- y^0}{h}s\right| - 1\right)}}{ \left |\frac{x- y^0}{r}-\frac{y- y^0}{h}s\right|},
\end{equation}
can be viewed as an easy-to-interpolate analytic function of
$(s, \theta, \varphi)$ around $s=0$ ($r=\infty$), where $\theta$ and
$\varphi$ are certain angular variables. In detail, using the
``singularity resolving'' spherical-coordinate change of variables
\begin{equation}\label{eq:defparametrizationins}
  {\mathbf x}(s, \theta, \varphi ;y^0) = 
  \tilde {\mathbf x}(h/s, \theta, \varphi ;y^0),
\end{equation}
where $\tilde {\mathbf x}$ denotes the $y^0$-centered classical
spherical-coordinate change of variables
\begin{equation} \label{eq:defparametrizationr}
\tilde {\mathbf x}(r, \theta, \phi ;y^0) = y^0 + \begin{pmatrix} r\sin \theta \cos \phi \\ r\sin \theta \sin \phi \\ r\cos \theta\end{pmatrix}, \qquad 0\leq r< \infty, \mkern5mu 0\leq \theta \leq \pi, \mkern5mu 0\leq \phi < 2\pi,
\end{equation}
the ``analytic factor'' $W_1(s,\theta,\phi,y;y^0)$, for
$y,y^0\in\mathbb{R}^3$, is defined by
\begin{equation}
  \label{eq:W_s_thet_phi}
  W_1(s,\theta,\phi,y;y^0) = \widehat{W}_1(\mathbf{x}(s, \theta, \varphi),y;y^0)
\end{equation}
so that,
\begin{equation}
  \label{eq:SL_thet_phi}
\mbox{for}\quad x =\mathbf{x}(s, \theta, \varphi ;y^0),   \quad\mbox{we have}\quad \Phi( x, y )= \frac{1}{4\pi} \frac{e^{\ii kr}}{r}W_1(s,\theta,\phi,y;y^0).
\end{equation}
Since $\left|\frac{x-y^0}{r} \right|=1$, $\frac{y- y^0}{h}\leq 1$, the
function $W_1$ is analytic for all $|s|<1$, including $s=0$
($r=\infty$). But lines~13 and~18 in Algorithm~1 (presented in
Section~\ref{sec:IFGF_rec_interp} below), show that the IFGF method
only interpolates onto points for which $r\geq \frac{3 H}{2}$, and,
thus, that the IFGF only uses values of $s$ satisfying
\begin{equation}\label{s_rest}
  s=\frac hr = \frac{\frac{\sqrt{3} H}{2}}{r}\leq \frac{\frac{\sqrt{3} H}{2}}{\frac{3 H}{2}}
  = \frac{\sqrt{3}}{3}\approx 0.577.
\end{equation}
In particular, the IFGF application utilizes the function $W_1$ well
within its domain of analyticity. This function is additionally slowly
oscillatory in the $s$ variable~\cite{BauingerBruno2021},
cf. Figure~\ref{fig:FactUnfaDL} and its caption.  The function $W_1$
may thus be efficiently interpolated in the variables
$(s,\theta,\phi)$, with a finite interpolation interval in the $s$
variable accounting for all functional variations up to an including
$r=\infty$, by using just a few interpolation intervals in the
variable $s$ in conjunction with adequately selected compact
interpolation intervals in the $\theta$ and $\phi$ variables, as
detailed in Section~\ref{sec:IFGF_interp}. As indicated in that
section, following~\cite{BauingerBruno2021}, the IFGF method relies on
this property to accelerate the evaluation of the
sum~\eqref{eq:field1} with kernel $K = \Phi$.

In what follows we seek an analogous factorization for the terms $V_2$
and $V_3$ in the double-layer expression
$\frac{\partial \Phi(x,y)}{\partial \nu(y)} = \frac{1}{4\pi}V_2(x,y) -
\frac{\ii k}{4\pi} V_3(x,y)$ displayed in
equation~\eqref{single_double_terms_2}. Using once again the notations
$r = |x- y^0|$ and $s=\frac{h}{r}$ we now obtain
\begin{subequations}\label{double_layer_fact}
\begin{align}
  V_2(x,y) & = \frac{e^{\ii kr}}{r^2}\frac{e^{\ii k\frac{h}{s}\left(\left |\frac{x-y^0}{r}-\frac{y-y^0}{h}s\right| - 1\right)}}{ \left |\frac{x-y^0}{r}-\frac{y-y^0}{h}s\right|^3} \left\langle \frac{x-y^0}{r}-\frac{y-y^0}{h}s, \nu(y) \right\rangle =\frac{e^{\ii kr}}{r^2}  \widehat{W}_2(x,y;y^0),\\
  V_3(x,y) & =  \frac{e^{\ii kr}}{r}\frac{e^{\ii k\frac{h}{s}\left(\left |\frac{x-y^0}{r}-\frac{y-y^0}{h}s\right| - 1\right)}}{ \left |\frac{x-y^0}{r}-\frac{y-y^0}{h}s\right|^2}\left\langle \frac{x-y^0}{r}-\frac{y-y^0}{h}s, \nu(y) \right\rangle = \frac{e^{\ii kr}}{r} \widehat{W}_3(x,y;y^0),
\end{align}
\end{subequations}
where
\begin{subequations}\label{eq:W2W2}
\begin{align}
  \widehat{W}_2(x,y;y^0) &= \frac{e^{\ii k\frac{h}{s}\left(\left |\frac{x-y^0}{r}-\frac{y-y^0}{h}s\right| - 1\right)}}{ \left |\frac{x-y^0}{r}-\frac{y-y^0}{h}s\right|^3} 
                                   \left\langle \frac{x-y^0}{r}-\frac{y-y^0}{h}s, \nu(y) \right\rangle \\
  \widehat{W}_3(x,y;y^0) &= \frac{e^{\ii k\frac{h}{s}\left(\left |\frac{x-y^0}{r}-\frac{y-y^0}{h}s\right| - 1\right)}}{ \left |\frac{x-y^0}{r}-\frac{y-y^0}{h}s\right|^2}
                                    \left\langle \frac{x-y^0}{r}-\frac{y-y^0}{h}s, \nu(y) \right\rangle.
\end{align}
\end{subequations}
Defining, for $y^0\in\mathbb{R}^3$ and $y\in B(y^0,H)$, the analytic
factors
\begin{equation}
  \label{eq:W_s_thet_phi_22}
  W_2(s,\theta,\phi,y;y^0) = \widehat{W}_2({\mathbf{x}(s, \theta, \varphi),y;y^0)}\quad\mbox{and}\quad W_3(s,\theta,\phi,y;y^0) = \widehat{W}_3({\mathbf{x}(s, \theta, \varphi),y;y^0)}
\end{equation}
the relations (cf.~\eqref{eq:W_s_thet_phi}),
\begin{equation}\label{eq:W_23_s_thet_phi}
\mbox{for}\quad x =\mathbf{x}(s, \theta, \varphi ;y^0),\quad  \frac{\partial \Phi(x,y)}{\partial \nu(y)} = \frac{1}{4 \pi}\frac{e^{\ii kr}}{r^2} W_2( s, \theta, \phi, y; y^0 ) 
                                       - \ii k \frac{1}{4 \pi}\frac{e^{\ii kr}}{r}   W_3( s, \theta, \phi, y; y^0 )
\end{equation}
follow. Clearly, like the analytic factor $W_1$ in
equation~\eqref{eq:W_s_thet_phi}, $W_2$ and $W_3$ are analytic
functions of $s$ for all $|s|<1$. Thus, per equation~\eqref{s_rest}
and associated parenthetical comment, these functions can be
effectively interpolated in the IFGF context using a small number of
interpolation intervals in the variable $s$ in conjunction with
adequately selected interpolation intervals in the $\theta$ and $\phi$
variables, as detailed in Section~\ref{sec:IFGF_interp}.

It is important to note here that, in many cases, it is more efficient
to evaluate double-layer potential contributions via interpolation of
the single function
\begin{equation}\label{w23}
  W_4 = \widetilde{W}_2- \ii k W_3,\quad\mbox{where}\quad  \widetilde{W}_2 =\widetilde{W}_2( s, \theta, \phi, y; y^0 ) = \frac{1}{r} W_2( s,
  \theta, \phi, y; y^0 ) =\frac{s}{h} W_2( s, \theta, \phi, y; y^0 ),
\end{equation}
instead of separately interpolating each of the two functions $W_2$
and $W_3$ (as required when using the
decomposition~\eqref{eq:W_23_s_thet_phi}), thereby leading to an IFGF
cost reduction by a factor of two.  Indeed, paralleling the
decompositions~\eqref{eq:W_s_thet_phi} and~\eqref{eq:W_23_s_thet_phi},
in this case we obtain the expression
\begin{equation}\label{eq:Wtilde_23_s_thet_phi}
  \frac{\partial \Phi(x,y)}{\partial \nu(y)} = \frac{1}{4 \pi}\frac{e^{\ii kr}}{r} \left(\widetilde{W}_2( s, \theta, \phi, y; y^0 )- \ii k    W_3( s, \theta, \phi, y; y^0 )\right) =  \frac{1}{4 \pi}\frac{e^{\ii kr}}{r} W_4( s, \theta, \phi, y; y^0 )
\end{equation}
in terms of the single function $W_4$. Like the previously considered
functions $W_1$, $W_2$ and $W_3$, the function $W_4$ is also an
analytic function of $s$ for $|s|<1$. As discussed in
Section~\ref{sec:CombLayIFGF}, interpolation of the function $W_4$ can
lead to larger errors in the resulting values of
$\frac{\partial \Phi(x,y)}{\partial \nu(y)}$, for small values of $h$,
but its use is advantageous for values of $h$ bounded away from zero.

The oscillatory character of the various analytic factors is
illustrated in Figure~\ref{fig:FactUnfaDL}---omitting, for brevity,
the factor $W_1$, which is considered in~\cite{BauingerBruno2021}, and
which is quite similar in form and character to $W_3$. For clarity,
for these examples the cube side $H$ is used as a measure of the size
of the interpolation problem instead of the parameter
$h = \frac{\sqrt{3}}{2}H$. The test setup is laid down at the top of
Figure~\ref{fig:FactUnfaDL}. It shows a single source point at the
position $y$ that, as it happens, gives rise to the fastest possible
oscillations along the measurement line (shown as a solid yellow line)
among all possible source positions at and within the blue source box.
Panels~(a) and~(b) demonstrate the highly oscillatory character of the
double-layer kernel components $V_2$ and $V_3$, even for small values
of $H$ ($H\leq 2\lambda$); the component $V_1$ exhibits a similar
oscillatory character. Panels (c)-(f), in turn, illustrate the
corresponding slowly oscillatory character of various IFGF analytic
factors introduced in the present section, even for much larger values
of $H$---e.g., for $H$ values up to $H=20\lambda$. Clearly, the
analytic factors $W_j$, $j=1,\dotsc,4$, oscillate much more slowly
than the corresponding full kernel components $V_j$, and, as expected
from the discussion above in this section, they actually tend to a
finite limit as $r\to\infty$.  All of these functions can therefore be
interpolated more effectively than the corresponding functions
$V_j$. As indicated above, however, a key difference in the analytic
factor $\widetilde{W}_2$ (and, thus, $W_4$), which is not noticeable
in this figure, becomes apparent as the problem of interpolation is
considered---as discussed in the following section and highlighted in
Figure~\ref{fig:IFGFSLandDL}.

\subsection{Factored Green function: single-box interpolation}\label{sec:IFGF_interp}
To appropriately structure an interpolation approach that takes
advantage of the slowly oscillatory character of the analytic factors
introduced in the previous section, we consider the $y^0$-centered
axes-aligned cubic box $B(y^0,H)$ of side $H$ and centered at $y^0$.
As discussed in what follows, for all target points $x_\ell$ that are
at least one box away from $B(y^0,H)$, the sum of all contributions
arising from the sources contained in $B(y^0,H)\cap\Bdy_N$ can be
effectively evaluated via interpolation of the corresponding analytic
factors $W_j$, $j=1,\dotsc,4$---since, for such target points we have
$(r,\theta,\phi) \in \mathcal{E} \coloneqq [3H/2,\infty] \times
[0,\pi] \times [0,2\pi)$ or, equivalently,
$(s,\theta,\phi) \in [0,\eta] \times [0,\pi] \times [0,2\pi)$, where
$\eta = \sqrt{3}/3$.

To tackle the interpolation task, the target domain $\mathcal{E}$ is
partitioned as described in what follows. Letting, for given positive
integers $n_s$ and $n_C$,
\begin{equation}\label{deltas}
  \Delta s = \frac{\eta}{n_s} \quad \text{ and } \quad 
  \Delta \theta = \Delta \phi = \frac{\pi}{n_C},
\end{equation}
and
$K_C \coloneqq \{1,\dotsc,n_s \} \times \{1,\dotsc,n_C\} \times
\{1,\dotsc,2n_C\}$, and letting, for each
$\bgamma = (\gamma_1,\gamma_2,\gamma_3)\in K_C$,
$s_{\gamma_1} = \gamma_1\Delta s$,
$\theta_{\gamma_2} = \gamma_2\Delta\theta$ and
$\phi_{\gamma_3} = \gamma_3\Delta\phi$, interpolation intervals
$E^s_{\gamma_1},\, E^{\theta}_{\gamma_2,\gamma_3}$ and
$E^{\phi}_{\gamma_3}$ along the $s,\,\theta,$ and $\phi$ directions
are defined by
\begin{equation}\label{eq:ConeIntervals}
  E^s_{\gamma_1} = [s_{\gamma_1-1},s_{\gamma_1}), \qquad
  E^{\theta}_{\gamma_2,\gamma_3} = 
  \begin{cases}
  [\theta_{n_C-1},\pi]                       & \text{for } \gamma_2 = n_C,\, \gamma_3 = 2n_C\\
  (0,\Delta\theta)                           & \text{for } \gamma_2 = 1,\,  \gamma_3 > 1\\
  [\theta_{\gamma_2-1},\theta_{\gamma_2})    & \text{otherwise,} 
  \end{cases} \qquad
  E^{\phi}_{\gamma_3} = [\phi_{\gamma_3-1},\phi_{\gamma_3}).
\end{equation}
For each $\bgamma = (\gamma_1,\gamma_2,\gamma_3) \in K_C$, we call
\begin{equation}\label{eq:ConeDomains}
  E_{\bgamma} \coloneqq E^s_{\gamma_1} \times E^{\theta}_{\gamma_2,\gamma_3} 
                 \times E^{\phi}_{\gamma_3}
\end{equation}
a \emph{cone domain}, and the image of $E_{\bgamma}$ under the
parametrization~\eqref{eq:defparametrizationins},
\begin{equation}\label{eq:ConeSegmentSet}
  C_{\bgamma}(y^0)\coloneqq \left\{ x = 
    \mathbf{x}(s,\theta,\phi ;y^0) \,:\, (s,\theta,\phi) \in E_{\bgamma} \right\},
\end{equation}
a \emph{cone segment}. Note that
\[
  \mathcal{E} = \bigcup \limits_{\bgamma \in K_C} C_{\bgamma} (y^0)
  \quad \text{ and } \quad
  C_{\bgamma}(y^0) \cap C_{\bgamma'}(y^0) = \emptyset \quad \text{ for } \bgamma \neq \bgamma'.
\]
\begin{remark}\label{cone_seg_size}
  The recursive interpolation method introduced in the following
  section, which utilizes boxes $B(y^0,H)$ of various sizes $H$ and
  for various centers $y^0$, requires, for each box size $H$ used,
  selection of numbers $n_s = n_s(H)$ and $n_C = n_C(H)$ in
  equation~\eqref{deltas} which ensure cone-segment interpolation at
  essentially constant accuracy for all used values of $H$. Such
  selection is facilitated by Theorems~1 and~2
  in~\cite{BauingerBruno2021}, which relate the acoustic size $kH$ of
  the box $B(y^0,H)$ to the error arising from cone-segment
  interpolation. In particular, these theorems imply that fixed
  numbers $n_s$ and $n_C$ give rise to essentially constant
  cone-segment interpolation errors for all $H$-side boxes satisfying
  $kH <1$. Additionally, these results tell us that, for box sizes in
  the complementary $kH\geq 1$ case, doubling the box size $H$
  requires doubling the numbers $n_s$ and $n_C$ (resulting in a
  two-fold refinement of cone segments in each of the three
  spherical-coordinate directions) in order to ensure an essentially
  unchanging cone-segment interpolation error. Accordingly, given
  integer initial values $n_{C,0}\geq 1$ and $n_{S,0}\geq 1$, the IFGF
  algorithm uses constant values $n_C=n_{C,0}$ and $n_s=n_{S,0}$ for
  all box sizes $H$ satisfying $kH\leq 1/2$, and, then, for box sizes
  $H$ for which $kH > 1/2$, the algorithm doubles the numbers $n_s$
  and $n_C$ as it transitions from a level $d$ (with, say, box size
  $H$) to a subsequent level $(d-1)$ (with box size $2H$). Throughout
  this paper the box sizes at the initial level $d=D$ were taken to
  satisfy the condition $kH< 1$, and the selections $n_{C,0}=2$ and
  $n_{S,0} = 1$ were used. Note that, for such values, at  levels
  $d$ for which $kH\leq 1/2$, if any such levels exist, each box
  $B(y^0,H)$ is furnished with a total of eight cone segments, each
  one an unbounded octant section ranging from $r= \frac{3H}{2}$
  ($s=\frac{\sqrt{3}}{3}$) to $r=\infty$ ($s=0$). Such a coarse
  cone-segment partition used in conjunction with polynomial
  interpolation by a single polynomial in each conical-segment
  coordinate variable, as described in what follows, with polynomial
  degrees such as, e.g., degree 3 in the radial variable $s$, and
  degree 5 in the angular variables $\theta$ and $\phi$, suffices to
  provide accuracies of the orders demonstrated in
  Section~\ref{sec:NumRes}.
\end{remark}

The proposed 3D interpolation strategy within cone segments
$C_{\bgamma}$ amounts to a tensor product of 1D Chebyshev
interpolation schemes based on use of fixed numbers $P_s$ and
$P_{\rm{ang}}$ of interpolation points along the radial and angular
variables, respectively, for a total of $P=P_s P^2_{\rm{ang}}$
interpolation points in each cone segment $C_{\bgamma}$. The
underlying 1D Chebyshev interpolating scheme is defined, for a given
function $u:[-1,1]\to \mathbb{C}$, as the interpolating polynomial
$I_n u$ of degree $(n-1)$ given by
\begin{equation}\label{eq:ChebyPoly}
  I_n u(x) = \sum_{i=0}^{n-1} a_i T_i(x), \qquad x \in [-1,1],
\end{equation}
where $T_i(x) = \cos(i\arccos(x))$. The coefficients $a_i$ are given by
\begin{equation}\label{eq:ChebyCoeffs}
  a_i = \frac{2}{b_i(n-1)} \sum_{k=0}^{n-1} c_k u(x_k) T_i(x_k), 
\end{equation}
where, for $i = 0,\dotsc,n-1$,
\begin{equation}\label{eq:ChebyNodes}
  x_k = \cos\left( \frac{2k+1}{2n} \pi \right), \quad 
  b_i = \begin{cases} 
        1 &  i\neq 0 \\
        2 &  i = 0, 
        \end{cases} \quad
  c_k = \begin{cases} 
        0.5 &  k = 0 \text{ or } k = n-1 \\
        1   &  \text{otherwise}.
        \end{cases}
\end{equation}

Use of the spherical-coordinate Chebyshev interpolation approach associated
with the box $B(y^0,H)$ can be applied to obtain approximate values of each one
of the analytic factors $W_j( s, \theta, \phi, y; y^0 )$, $j=1,\dotsc,4$
defined in Section~\ref{sec:GenIFGF}, for any given source point $y\in
B(y^0,H)$ and any given target point $x = \mathbf{x}( s, \theta, \phi)$ (with
$y^0$-centered spherical coordinates $( s, \theta, \phi)$) that is at least one
size-$H$ box away from $B(y^0,H)$.  By linearity, the same interpolation method
can be applied, with similar accuracy, to linear combinations such as the one
given by the sum over $y$ on the right-hand side of the relation
\begin{equation}\label{exponential_sum_linearity}
  S^{H,y^0}_j(s,\theta,\phi) = \frac{e^{\ii k |x-y^0|}}{4\pi |x-y^0|^{q_j}}
   \sum\limits_{y \in B(y^0,H) \cap \Bdy_N}  a(y)W_j(s,\theta,\phi,y;y^0),
\end{equation}
for any given number of sources within the box $B(y^0,H)$ and for
given coefficients $a(y)$. Here, with reference to the factorizations
introduced in Section~\ref{sec:GenIFGF}, we have set $q_2 = 2$ and
$q_j = 1$ for $j=1,3,4$.

In view of the factorizations displayed in~\eqref{eq:SL_thet_phi}
and~\eqref{eq:W_23_s_thet_phi}, interpolation of the sum over $y$
in~\eqref{exponential_sum_linearity} with $j=1,2,3$ can be used to
accelerate the evaluation of the sum~\eqref{eq:field1} with $K$ equal
to either the single-layer or the double-layer
kernels~\eqref{eq:GreenFunction} and~\eqref{double_layer}. As
discussed in the following section, however, in the case of the
double-layer kernel a more efficient interpolation strategy can be
devised by utilizing, in addition, the
factorization~\eqref{eq:Wtilde_23_s_thet_phi} associated with the
analytic factor $W_4$.

\subsection{IFGF method for the double-layer
  potential}\label{sec:CombLayIFGF}

\begin{figure}[h!]
    \centering
    \includegraphics[width=0.8\textwidth]{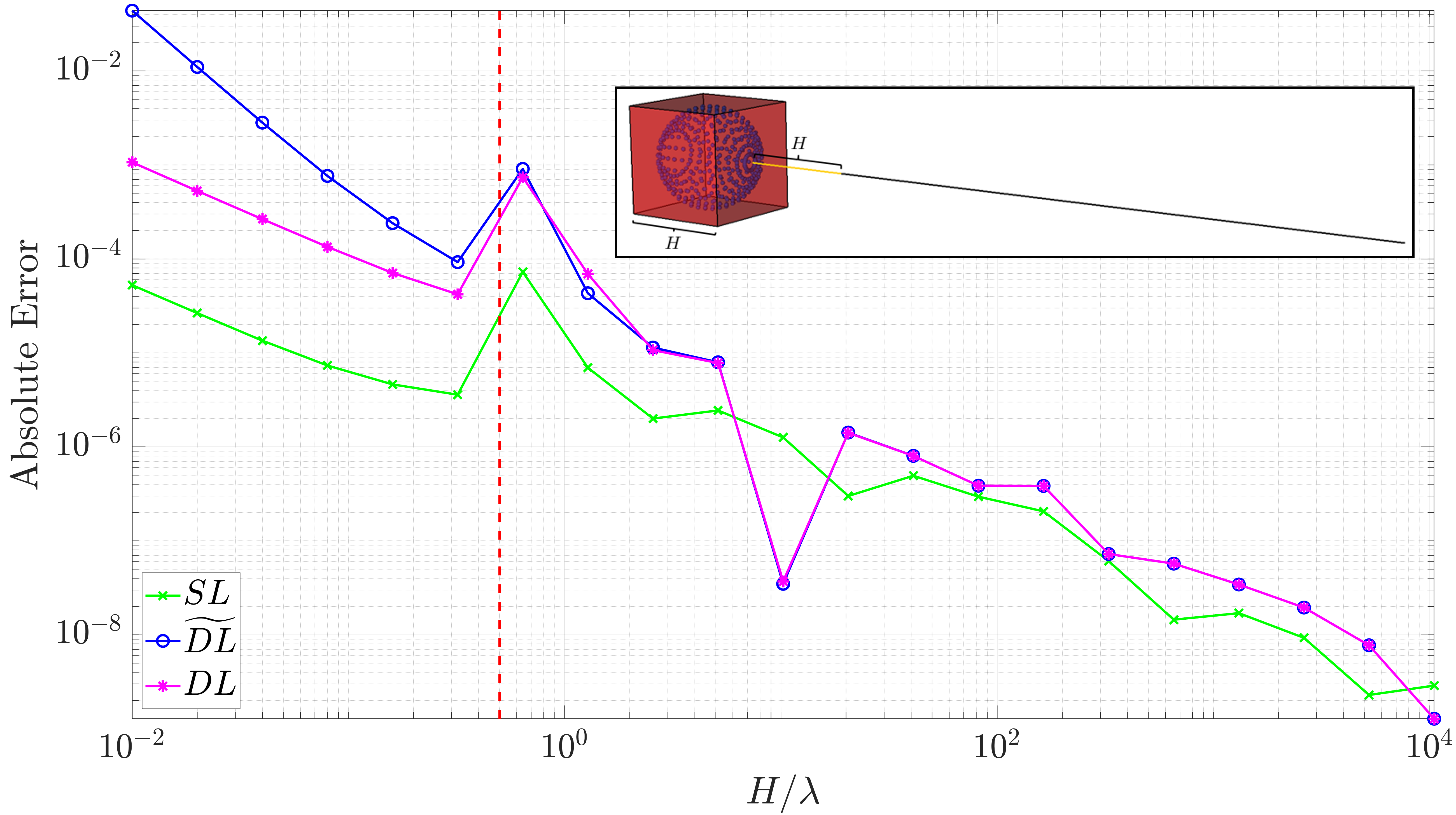}
    \caption{\small Maximum absolute errors arising in IFGF
      representation of the single-layer and double-layer
      potentials. In the configuration considered for this test, 200
      sources $y \in \Gamma\subset B(y^0,H)$, located on the spherical
      surface $\Gamma$ and contained in the box $B(y^0,H)$, were used,
      as shown in the inset. The curves presented illustrate the
      character of three IFGF-interpolation algorithms, namely,
      1)~Interpolation of the single-layer potential on the basis
      of~\eqref{eq:SL_thet_phi} by means
      of~\eqref{exponential_sum_linearity} with $j=1$ (curve labeled
      ``$SL$'', shown in green); 2)~Interpolation of the double-layer
      potential on the basis of
      equation~\eqref{eq:Wtilde_23_s_thet_phi} by means
      of~\eqref{exponential_sum_linearity} with $j=4$ (curve labeled
      ``$\widetilde{DL}$'', shown in blue), and 3)~Interpolation of
      the double-layer potential on the basis of
      equation~\eqref{eq:W_23_s_thet_phi} by means of a combination of
      the $j=2$ and $j=3$ versions
      of~\eqref{exponential_sum_linearity} (curve labeled ``$DL$'',
      shown in magenta).  The errors shown were evaluated as the
      maximum of the errors obtained at a set of 200 target points
      $x=\mathbf{x}(s, \theta, \varphi; y^0)$, with
      $\theta=\varphi =0$, located at a distance greater than or equal
      to $H$ from the boundary of $B$ along a radial line exterior to
      $B$, as shown in the inset.  Clearly the $\widetilde{DL}$ and
      $DL$ errors essentially coincide to the right of the line
      $H/\lambda = 0.5$ (vertical dashed red line)---, which, as
      discussed in Section~\ref{sec:CombLayIFGF}, enables certain
      efficiency gains in the evaluation of double-layer potentials.}
  \label{fig:IFGFSLandDL}
\end{figure}

In order to study the relative advantages offered by the interpolation
strategies arising from the concepts introduced in
Sections~\ref{sec:GenIFGF} and~\ref{sec:IFGF_interp}, we consider
Figure~\ref{fig:IFGFSLandDL}, which displays absolute errors in IFGF
approximations associated with the single- and double-layer
potentials. For this experiment, single- and double-layer potentials
resulting from a total of 200 source points on a spherical surface
contained within the box $B(y^0,H)$ were considered, as depicted in
the figure inset, and errors along the radial line shown in the inset
were evaluated and are displayed in the figure. Only radial
interpolations were performed, since the angular factorizations are
identical for all the factored forms used. The various error curves,
labeled \textit{SL}, \textit{DL} and $\widetilde{DL}$ in the figure,
display absolute errors arising in the approximation of the single-
and double-layer potentials, namely, the error \textit{SL} in the
single-layer approximation resulting from use of the
factorization~\eqref{eq:SL_thet_phi}, and the errors \textit{DL} and
$\widetilde{DL}$ in the double-layer approximations resulting from use
of the factorizations~\eqref{eq:W_23_s_thet_phi}
and~\eqref{eq:Wtilde_23_s_thet_phi}, respectively. Mirroring the
octree IFGF structure described in the following section, abscissae
for the points in the graph are given by powers of $2$ times a fixed
constant. In each case the error resulting from interpolation of
analytic factors was computed over a single interpolation interval,
namely, the interpolation interval starting at a distance $H$ from the
source box---the closest interval to the source box, and, therefore,
the one giving rise to the largest interpolation errors. In accordance
with Remark~\ref{cone_seg_size}, for $kH < 1$, the interval spans the
complete semi-infinite line (i.e., the complete interval
$0 \leq s \leq\frac{\sqrt{3}}{3})$), and, starting from $k H = \pi$
(that is, $H/\lambda = 0.5$, which is marked by a vertical dashed line
in the figure) the interval size in the $s$ variable is halved every
time $kH$ is doubled.  Interpolations were performed using a total of
$P_s=3$ interpolation points along the interval, and IFGF errors were
evaluated at $200$ equispaced target test points in the
interval. Clearly, over wide $H/\lambda$ ranges, spanning from
$H/\lambda = 0.5$ all the way to $H/\lambda =\infty$, the
interpolation scheme associated with the factorization
$\widetilde{DL}$, which requires interpolation of the single analytic
factor $W_4$, incurs the same error as the more expensive scheme
arising from interpolation of the decomposition \textit{DL}---which
requires interpolation of the two analytic factors $W_2$ and
$W_3$. The curve indicates, however, that the $\widetilde{DL}$ scheme
does lead to significant accuracy losses versus the \textit{DL}
scheme, when used for $H/\lambda < 0.5$.  The increased error
displayed by the $\widetilde{DL}$ curve in the range $H/\lambda < 0.5$
can be traced to the partial factorization of the $\frac{1}{r}$ term
associated with the analytic factor $\widetilde{W}_2$ (which still
contains a singular component, as evidenced by the $1/r$ term in
equation~\eqref{w23}), versus the fully regularized
decomposition~\eqref{eq:W_23_s_thet_phi}. The right-hand expression in
equation~\eqref{w23} provides a different perspective on the
phenomenology, as it shows that the $\widetilde{DL}$ factorization
suffers from increased errors for small values of $h$, and thus, in
view of~\eqref{eq:def_eps}, for small values of $H/\lambda$.

While rigorously accounting for the Green function singularity and
maintaining accuracy for arbitrarily small values of the source box
size $H/\lambda$, the approach demonstrated in the \textit{DL} curves,
which is based on interpolation of the analytic factors $W_2$ and
$W_3$, doubles the cost of the IFGF accelerator vis-\`a-vis the
$W_4$-based method associated with the $\widetilde{DL}$ curves. An
optimal strategy then calls for use of the $W_2$-$W_3$ based method
for the portion of the IFGF octree containing boxes of size
$H/\lambda < 0.5\lambda$, and use of the $W_4$-based method for the
remainder the octree.  For the structures considered in this paper,
most of which do not contain significant subwavelength geometric
features, the second approach is advantageous throughout the complete
IFGF octree---as it enjoys the reduced computational cost while
preserving accuracy. In a more general context a hybrid
$W_2$-$W_3$/$W_4$ approach as described above would be used which,
additionally, would incorporate adaptivity in the box sizes in regions
containing finely meshed surfaces (as discussed in
Section~\ref{sub_scat} in connection with the submarine propeller
structure). Such extensions are beyond the scope of this paper, and
are left for future work.

%

\subsection{Factored Green function: recursive
  interpolation}\label{sec:IFGF_rec_interp}
\begin{figure}
    \centering
    \includegraphics[width=\textwidth]{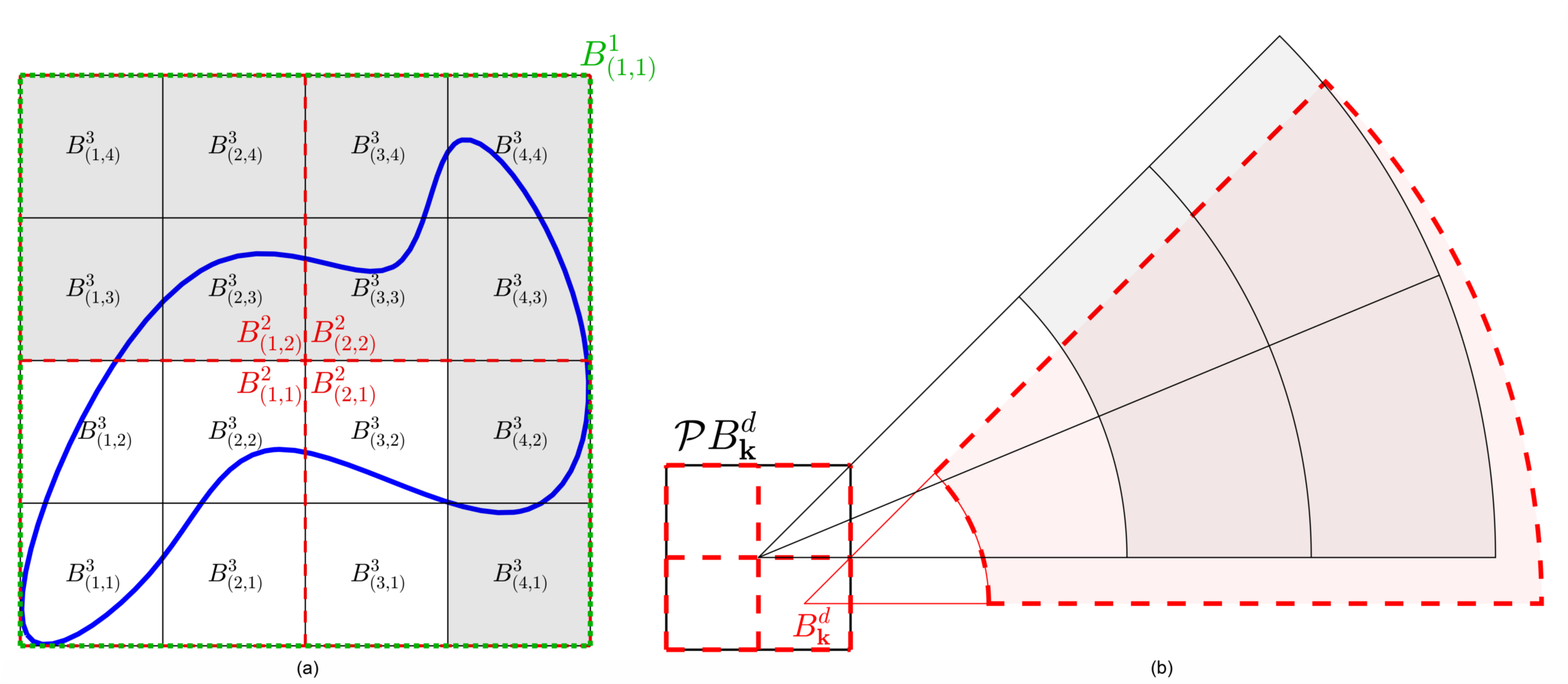}
    \caption{\small (a) Two-dimensional sketch of the IFGF domain
      decomposition for the three-level $(D=3)$ case, with a scatterer
      is sketched in blue. In particular, the figure displays the
      neighbors (white boxes) and cousins (gray boxes) of the box
      $B_{(2, 1)}^3$. (b) Illustration of typical cone segments for a
      box $B_\mathbf{k}^d$ (shown in red and centered at the southeast
      red box) and its parent box with some associated cone segments
      (drawn in black).}
    \label{fig:domaindecomposition}
\end{figure}

To effect the desired evaluation of discrete
operators~\eqref{eq:field1} at $\mathcal{O}(N\log N)$ computing cost,
the IFGF approach implements the interpolation scheme described in
Sections~\ref{sec:GenIFGF} through~\ref{sec:CombLayIFGF} in a
recursive manner, in sets of larger and larger boxes, so that the
values required by the interpolation scheme associated with a large
box can themselves be obtained by interpolation from smaller
boxes. The IFGF method implements such a recursive interpolation
scheme on the basis of a $D$-level octree sequence of partitions of
the discrete surface $\Bdy_N$~\cite{BauingerBruno2021} for a suitably
selected $D \in \N$.  Thus, the boxes are defined iteratively starting
from the single level $d=1$ box $B_{(1, 1, 1)}^1 \supset \Bdy_N$,
leading to the trivial partition of $\Bdy_N$ into a single subset
equal to all of $\Bdy_N$. For $d = 2, \dotsc, D$, in turn, the
level-$d$ boxes are defined via partition of each of the level $(d-1)$
boxes into eight equi-sized and disjoint boxes of side
$H_d = H_{d-1}/2$, resulting in the level $d$ boxes
$\Bkd = B(\ykd,H_d)\subset \mathbb{R}^3$ (see
equation~\eqref{eq:box}), for
$\mathbf{k} \in I^d_B \coloneqq \{1, \ldots, 2^{d-1}\}^3$, and for
certain box centers $\ykd$. Clearly, each box $\Bkd$ on level $d$
($2 \leq d \leq D$) is contained in a parent box on level $d-1$, which
we denote by $\mathcal{P}\Bkd$. The level-$d$ partition of $\Bdy_N$ is
obtained via intersection of $\Bdy_N$ with each one of the level-$d$
boxes.  Figure~\ref{fig:domaindecomposition}(a) illustrates an
analogous two-dimensional hierarchical quadtree structure in the
$D = 3$ (three-level) case.

On the basis of these concepts we utilize the functions
$S^{H,y^0}_j(s,\theta,\phi)$ in~\eqref{exponential_sum_linearity} with
$y^0 = \ykd$ and $H=H_d$ (and, thus, with $B(y^0,H) = \Bkd$), which in
what follows are
denoted by
\begin{equation}\label{exponential_sum_1}
  S^{\mathbf{k},d}_j(s,\theta,\phi) = S^{ H_d,y^{\mathbf{k}}_d}_j(s,\theta,\phi) 
  = \frac{e^{\ii k r_\mathbf{k}}}{4\pi r_\mathbf{k}^{q_j}}
  F^{\mathbf{k},d}_j(s,\theta,\phi),
\end{equation}
where $a(y)$ are given coefficients, where
$r_{\mathbf{k}} \coloneqq |\mathbf{x}(s, \theta, \varphi ;\ykd) -
y^{\mathbf{k}}_d|$, and where,
\begin{equation}\label{exponential_sum_11}
\mbox{letting}\quad \Wkd =
W_j(s,\theta,\phi,y;y^\mathbf{k}_d),\quad\mbox{we have set}\quad
F^{\mathbf{k},d}_j(s,\theta,\phi)= \sum\limits_{y \in B_\mathbf{k}^d
  \cap \Bdy_N} a(y)\Wkd.
\end{equation}
As discussed below in this section, the proposed recursive
interpolation approach relies on application of the single-box
interpolation strategy described in Sections~\ref{sec:GenIFGF}
through~\ref{sec:CombLayIFGF} to each one the source boxes
$\Bkd = B(y^{\mathbf{k}}_d,H_d)$ and factors $F^{\mathbf{k},d}_j$,
starting at level $d=D$ and then proceeding to levels $d=D-1,...$
until all relevant levels $d=D,\dotsc,3$ have been tackled. (The
algorithm indeed stops at level $d= 3$ since, per construction of the
boxes and definition of cousins, at that stage all boxes are cousins
and therefore, at the end of that stage all interactions have already
been taken into account.)

Clearly, it is only necessary to effect interpolations for
\textit{relevant boxes} $B_{\mathbf{k}}^d$, namely, boxes
$B_{\mathbf{k}}^d$ that contain at least one surface discretization
point; the set $\mathcal{R}_B$ of all relevant boxes, including boxes
at any level $d$ ($1\leq d\leq D$), is denoted by
\begin{equation}\label{eq:RelevantBoxes}
    \mathcal{R}_B := \left\{ \Bkd \, : \, \Bdy_N \cap B_{\mathbf{k}}^d 
    \neq \emptyset,\, 1 \leq d \leq D,\, \mathbf{k} \in I^d_B \right\}.
\end{equation}
To facilitate the discussion, several additional but related concepts are
introduced, including, for a given level-$d$ box $\Bkd$, the set $\mathcal{N}
\Bkd$ of all level-$d$ boxes that are \textit{neighbors} of $\Bkd$, and the set
$\mathcal{M} \Bkd$ of all level-$d$ boxes that are \textit{cousins} of $\Bkd$
(that is to say, non-neighboring boxes who are children of the parents'
neighbors). The sets $\mathcal{U} \Bkd$ and $\mathcal{V} \Bkd$ of
\textit{neighbor points} and \textit{cousin points} of $\Bkd$, respectively,
are defined as the set all points in $\Bdy_N$ that are contained in neighbor
and cousin boxes of $\Bkd$, respectively. We clearly have
\begin{subequations}\label{eq:NeighborsCousins}
\begin{align}
    \mathcal{N} \Bkd &:= \{ B_\mathbf{j}^d \in \mathcal{R}_B : ||\mathbf{j} - \mathbf{k}||_\infty \leq 1 \}, \\
    \mathcal{M} \Bkd &:= \{ B_\mathbf{j}^d \in \mathcal{R}_B : B_\mathbf{j}^d \notin \mathcal{N} B_\mathbf{k}^d \wedge \mathcal{P}B_{\mathbf{j}}^d \in \mathcal{N} \mathcal{P} B_\mathbf{k}^d \} \\
    \mathcal{U} \Bkd &:= \left( \bigcup \limits_{B \in \mathcal{N} B_\mathbf{k}^d} B \right) \cap \Bdy_N  \\
    \mathcal{V} \Bkd &:= \left( \bigcup \limits_{B \in \mathcal{M} B_\mathbf{k}^d} B \right) \cap \Bdy_N.
\end{align}
\end{subequations}
Figure \ref{fig:domaindecomposition}(a) presents, in a two-dimensional
setting, the neighbor and cousin boxes (shown in white and gray,
respectively) of the box $B_{(2, 1)}^3$.

The IFGF algorithm accelerates the evaluation of \eqref{eq:field1} by
incrementally incorporating pairwise interactions of level-$d$ cousin
boxes, for $d = D, \ldots, 3$: at stage $d$ the algorithm evaluates,
by cone-segment spherical-coordinate interpolation, contributions from
each level-$d$ box $\Bkd$ to all surface discretization points that
lie on cousin boxes of $\Bkd$, as well as all cone-segment
interpolation points necessary at the subsequent level $(d-1)$ of the
interpolation scheme.  Note that
when level $d$ is reached, values of fields resulting from sources
within $\Bkd$ at points in boxes neighboring $\Bkd$ have already been
computed at the previously treated levels $D,D-1,\dotsc,d+1$ (but see
Remark~\ref{level_D_neighboring} below).  Once the level-$d$
interpolations have been performed at all necessary cousin points and
level-$(d-1)$ interpolation points associated with the box $\Bkd$, the
complete result $S^{\mathbf{k},d}_j(s,\theta,\phi)$ for that box at
each required cousin and interpolation point is obtained by
multiplication of the interpolated values
$F^{\mathbf{k},d}_j(s,\theta,\phi)$ with the corresponding value of
the ``centered factor''
\begin{equation}\label{cent_fact}
E^{\mathbf{k},d}_j(s,\theta,\phi) =   \frac{e^{\ii k r_\mathbf{k}}}{4\pi r_\mathbf{k}^{q_j}}.
\end{equation}
\begin{remark}\label{level_D_neighboring}
  An exception to the interpolation-based method just described occurs
  at the finest (initial) level $d = D$, for which there are no
  preceding levels, and for which fields at boxes neighboring
  $B_{\mathbf{k}}^{D}$ have therefore not been previously
  obtained. Thus, at level $d=D$, interactions from source points
  $B_{\mathbf{k}}^{D}$ to points in neighboring boxes are not
  evaluated by means of an interpolation scheme. Instead, in the
  context of this section, in which the values of the
  sum~\eqref{eq:field1} are to be produced, such interactions are
  obtained via direct summation of the corresponding
  $B_{\mathbf{k}}^{D}$ contributions. In other contexts, in which the
  IFGF is used as part of an integration scheme (such as the one
  introduced in Section~\ref{sec:RecPolar}), certain nearby
  interactions need to be handled by specialized quadrature rules
  (such as the ones introduced in Section~\ref{sec:SingInter}) while
  others do proceed via direct summation.
\end{remark}

As suggested above, the IFGF method evaluates the necessary cousin
interactions for a single- and double-layer kernels by means of simple
piecewise interpolation of the relevant factor(s)
$F^{\mathbf{k},d}_j(s,\theta,\phi)$ ($j=1,\dots,4$). Per the
discussion in Sections~\ref{sec:GenIFGF}
through~\ref{sec:CombLayIFGF}, the interpolations associated with a
level-$d$ relevant box $\Bkd$ ($\mathbf{k}\in I_B^d $) are effected by
relying on level-$d$ cone segments, which in what follows are denoted
by
\[
  C_{\mathbf{k}; \bgamma}^d = C_{\bgamma}(\ykd),
  \quad \bgamma\in K_C^d, \, \mathbf{k} \in I^d_B
\]
(cf.~\eqref{eq:ConeSegmentSet}), where, for certain selections of
integers $n_s = n^d_s$ and $n_C = n^d_C$, we have used the index set
$K_C^d = \{1,\dotsc,n^d_s \} \times \{1,\dotsc,n^d_C\} \times
\{1,\dotsc,2n^d_C\}$, resulting in radial and angular interpolation
intervals of lengths $\Delta s = \frac{\eta}{n^d_s}$ and
$\Delta \theta = \Delta \phi = \frac{\pi}{n^d_C}$
(cf. equation~\eqref{deltas}). In practice, the integers $n^d_s$ and
$n^d_C$ are selected so as to ensure constant accuracy across levels
per the prescriptions outlined in Remark~\ref{cone_seg_size}---which,
as the algorithm progresses from level $d$ to level $(d-1)$, requires
(resp. does not require) refinements of cone-segments provided
$kH_d > 1/2$ (resp.  $kH_d\leq 1/2$).  For each box $\Bkd$ and each
$\bgamma\in K_C^d$, we denote by
$\mathcal{X} C_{\mathbf{k}; \bgamma}^d$ the set of interpolation
points within the cone segment $C_{\mathbf{k}; \bgamma}^d$.  As in the
case of boxes, wherein only relevant boxes are used in the
interpolation process, for each box $\Bkd$ the method only utilizes
the co-centered cone segments that belong to the set
$\mathcal{R}_C B_\mathbf{k}^d$ of \textit{relevant cone segments},
defined by
\begin{subequations}\label{eq:RelevantCones}
\begin{align}
    \mathcal{R}_C B_\mathbf{k}^d &:= \emptyset \quad \text{for} \quad d \in \{1, 2\}, \\
    \mathcal{R}_C B_\mathbf{k}^d &:= \left\{ C_{\mathbf{k}; \bgamma}^d \, : \, \nu \in I_C^d\, , \, C_{\mathbf{k}; \bgamma}^d \cap \mathcal{V} B_\mathbf{k}^{d} \neq \emptyset \text{ or } C_{\mathbf{k}; \bgamma}^d \cap \left( \bigcup \limits_{C \in \mathcal{R}_C \mathcal{P} B_\mathbf{k}^{d}} \mathcal{X} C \right) \neq \emptyset \right \} \quad \text{for} \quad d \geq 3.
\end{align}
\end{subequations}
A two-dimensional sketch of two box-centered cone segments for a box
$B_\mathbf{k}^d$ and its parent $\mathcal{P}B_\mathbf{k}^d$ is
presented in Figure \ref{fig:domaindecomposition}(b). It is useful to
note that, for a given surface $\Bdy_N$, wavenumber $k$, and
prescribed error tolerance, the relevant cone segments can be
determined recursively---but in reverse order, starting from $d=3$ and
moving upwards the tree to $d=D$~\cite{BauingerBruno2021}.
\begin{remark}\label{cone_error}
  Following the cone-segment interpolation prescriptions described in
  Remark~\ref{cone_seg_size}, a given interpolation error tolerance
  will be met at every level $d\leq D$ provided it is met at level
  $D$---which can be ensured via appropriate selections of the
  cone-partitioning numbers $n_s$ and $n_C$ and the number of
  interpolation points $P_s$ and $P_{\rm{ang}}$ used.  Error
  accumulation over the $D = \log_2N$ stages of the method needs also
  to be taken into account to ensure an overall error tolerance is
  satisfied. In practice, setting the number of interpolation points
  for the $s$ and angular variables to $P_s = 3$ and
  $P_{\rm{ang}} = 5$, respectively, are sufficient to ensure an
  overall accuracy of the order of $10^{-4}$ for problems up to
  $D= 10$ levels and beyond.
\end{remark}
\begin{figure}
    \centering
    \includegraphics[width=\textwidth]{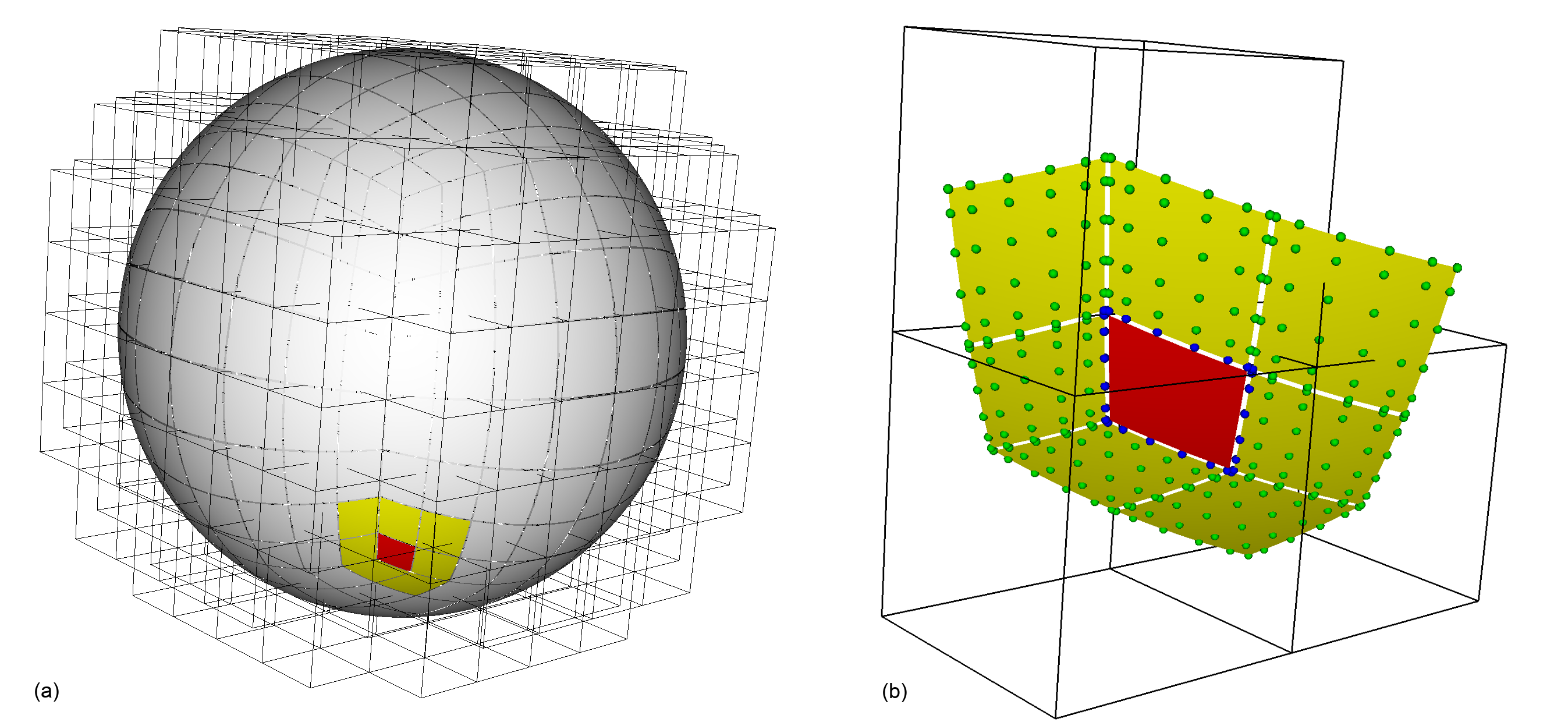}
    \caption{\small (a) Four level $(D=4)$ IFGF domain decomposition
      for a sphere (only the finest level boxes are shown). A source
      surface patch and its neighboring patches are highlighted in red
      and yellow, respectively. (b) Near-singular and regular points
      neighboring the red patch are depicted in blue and green,
      respectively.  All interactions from the red patch to itself are
      singular and, together with the near-singular and neighboring
      regular interactions are computed by means of the
      rectangular-polar algorithm described in
      Section~\ref{sec:RecPolar}. The remaining interactions are
      accelerated by means of the IFGF algorithm.}
    \label{fig:LocalInter}
\end{figure}

To achieve fast computation times, the IFGF method relies on the
factorization~\eqref{exponential_sum_1}, and it proceeds, as indicated
above, by interpolation of the factor
$F^{\mathbf{k},d}_j(s,\theta,\phi)$.  In view of the analyticity and
the slowly-oscillating character of the analytic factors
$W^{\mathbf{k},d}_j$ in~\eqref{exponential_sum_11}, accurate
interpolation of the factors $F^{\mathbf{k},d}_j(s,\theta,\phi)$ can
be achieved by polynomial interpolation in relatively coarse
$(s,\theta,\phi)$ meshes. In keeping with~\cite{BauingerBruno2021},
Chebyshev polynomials of degrees $P_s = 3$ and $P_{\rm{ang}} = 5$
points are used in this paper for interpolation along the radial and
angular directions, respectively, but, as indicated in that reference,
interpolation methods of higher order may of course be employed as
well.

For a given point $x$, the IFGF algorithm produces the desired field
$I_K(x)$ in \eqref{eq:field1} by subsequently incorporating, at each
level $D, D-1, D-2,\dotsc,3$, all relevant contributions from boxes
for which $x$ is a cousin point. (As indicated above, the
contributions from nearest neighbor boxes in level $D$ are computed
without recourse to the IFGF approach, using directly the methods in
Section~\ref{sec:RecPolar}.) Contributions from boxes for which $x$ is
not a cousin point at a given level are incorporated at subsequent
levels by grouping of boxes into larger boxes. The data necessary for
interpolation from a given box is obtained by adding all contributions
from the relevant children of that box---each one of which is
obtained, once again, by interpolation.

The details of the IFGF algorithm are presented in Algorithm
\ref{alg:ifgf}. Note that the algorithm does not compute the ``local
interactions'', that is, the interactions between neighboring boxes on
the finest level $(d=D)$. For the scattering solver proposed in this
paper, such interactions are evaluated by means of two separate
methods, as illustrated in Figure~\ref{fig:LocalInter} and described
in what follows. If the distance from a box's source points to a
neighboring box's target point $x_{\ell}$ is larger than the proximity
distance $\delta$ defined in Section~\ref{sec:RecPolar}, the local
interactions are evaluated by means of the algorithm described in
Section~\ref{sec:ReguInter}. In case the source points and the target
point in neighboring boxes are a distance greater than $\delta$, in
turn, the quadrature presented in Section~\ref{sec:SingInter} is used
instead. We emphasize that local interactions are computed \emph{only
  in the finest level} $(d=D)$.  Figure~\ref{fig:LocalInter}(a) shows
a four-level $(D=4)$ IFGF domain decomposition for a sphere, where a
source patch and its neighbors patches are highlighted in red and
yellow, respectively. A close-up view of these patches is shown in
Figure~\ref{fig:LocalInter}(b); source-to-neighbor near-singular
evaluation points are drawn in blue, while regular evaluation points
are depicted in green. All evaluations from a source patch to itself
are considered to be singular and are also computed using the
algorithms presented in Section~\ref{sec:SingInter}.
\begin{algorithm}[htp!] 
  \begin{algorithmic}[1]
    \small
    \State \textbackslash \textbackslash Direct evaluations of the
    highest-level local interactions.
    \For{$B_\mathbf{k}^D \in
      \mathcal{R}_B$} \label{algstate:looprelevantboxeslevelD}
    \For{$C_{\mathbf{k};\bgamma}^D \in \mathcal{R}_C
      B_\mathbf{k}^D$} \label{algstate:looprelevantconeslevelD}
    \Comment{Evaluate $F$ at all relevant interpolation points}
    \For{$x = \mathbf{x}(s, \theta, \varphi) \in \mathcal{X} 
      C_{\mathbf{k};\bgamma}^D$} \label{algstate:loopinterpointslevelD}
    \State Evaluate and store $F_\mathbf{k}^D(x)$. \EndFor \EndFor
    \EndFor
  \State
  \State \textbackslash \textbackslash Interpolation onto surface discretization points and parent
  interpolation points.
  \For{$d = D, \ldots, 3$}\label{algstate:loopd}
  \For{$B_\mathbf{k}^d \in \mathcal{R}_B$} \label{algstate:looprelevantboxes} 
  \For{$x \in \mathcal{V} B_\mathbf{k}^d$} \label{algstate:loopnearestneighbouringsurfacepoints}
  \Comment{Interpolate at cousin surface points} 
  \State Evaluate $I_\mathbf{k}^d(x)$ by interpolation \label{algstate:interpolationtosurfacepoints}
  \EndFor 
  \If {$d > 3$}
  \Comment{Evaluate $F$ on parent interpolation points} \State Determine parent $B_\mathbf{j}^{d-1} = \mathcal{P} B_\mathbf{k}^d$
  \For{$C_{\mathbf{j}; \bgamma}^{d-1} \in \mathcal{R}_C B_\mathbf{j}^{d-1}$} \label{algstate:looprelevantcones} 
  \For{$x \in \mathcal{X} C_{\mathbf{j};\bgamma}^{d-1}$} \label{algstate:loopinterppoints}
  \State Evaluate and add $F_\mathbf{k}^d(x) K(x, y_\mathbf{k}^d)/ K(x, y_\mathbf{j}^{d-1})$
  \EndFor 
  \EndFor 
  \EndIf 
  \EndFor 
  \EndFor
  \caption{\small IFGF Algorithm}
\label{alg:ifgf}
\end{algorithmic}
\end{algorithm}

\subsection{IFGF OpenMP parallelization}\label{sec:OpenMP}
The proposed OpenMP parallel IFGF
method~\cite{BauingerBruno_parallel_2021} is based in a fundamental
manner on an ordering of the IFGF box-cone structure. In particular,
the set
\begin{equation} \label{eq:relconesegmentslevelD} \mathcal{R}_C^d :=
  \bigcup \limits_{\mathbf{k} \in I_B^d: B_\mathbf{k}^d \in
    \mathcal{R}_B} \mathcal{R}_C B_\mathbf{k}^d, \quad \text{for } 3
  \leq d \leq D.
\end{equation}
(cf. Section~\ref{sec:IFGF_rec_interp}) of all level-$d$ cone segments
is ordered in accordance to a Morton order (which starts by ordering
the set of relevant boxes by means of a three-dimensional Z-looking
curve that zig-zags through all boxes~\cite{1993HashedOcttreeWarren},
and then orders the cone segments associated with each box along
angular and radial directions), so that the complete set of
interpolating polynomials at any given level $d$ is stored as a single
vector, sections of which are subsequently distributed to threads in
the OpenMP implementation. The main advantage provided by the IFGF
method is that the number of IFGF tasks at each level $d$, which
coincides with the number of elements in the set $\mathcal{R}_C^d$ of
level-$d$ cone segments, {\em is a large number that is essentially
  independent of $d$}---so that the set of tasks to be performed at
any given level $d$ remains {\em naturally} parallelizable for all
$d$. In other approaches such as the FMM, the acceleration strategy
requires use of decreasing number of increasingly larger structures
(as the octree is traversed from the leaves to the root) each one of
which is difficult to parallelize---leading to an associated {\em
  Parallelization
  Bottleneck}~\cite{2007Volakis,2007DirectionalFMMLexing,2014ParallelDirectionalFMMLexing}.
The IFGF approach, in contrast, incorporates large level-independent
numbers of OpenMP threads leading to a coarse-grained and
well-balanced parallelization resulting in a effective parallel
scaling at low memory cost. An extended discussion concerning
parallelization of the IFGF method, both under OpenMP and MPI
interfaces, can be found in~\cite{BauingerBruno_parallel_2021}.

\section{IFGF/RP Evaluation of surface integral
  operators} \label{sec:RecPolar}
The IFGF algorithm introduced in the previous section can be utilized
to accelerate the numerical evaluation of integral operators of the
form~\eqref{eq:gen-int-op} for any given discretization of the
scattering boundary $\Gamma$. The RP discretizations utilized in this
paper, in turn, are based on use of Chebyshev expansions over each one
of the patches $\Gamma^q$ introduced in Section~\ref{sec:SurfRep}. As
indicated in that section, the design of suitable quadrature methods
for evaluation of $\opI^q[\varphi](x)$ requires careful consideration
of the relative position of the point $x\in\Gamma$ and the integration
patch $\Gamma^q$. The ``regular'' case where $x$ is ``far'' from
$\Gamma^q$ is treated easily on the basis of regular Clenshaw-Curtis
quadrature (Section~\ref{sec:ReguInter} below) and is subsequently
accelerated by means of IFGF; the singular and near-singular cases,
for which either $x\in \Gamma^q$ or $x$ is ``near'' $\Gamma^q$, in
turn, are considered in Section~\ref{sec:SingInter}.

In detail, we consider a general integral operator $\opI^q$ with
singular kernel $K^q$ and density $\varphi^q$ defined over a component
patch $\Gamma^q$ is expressed in the parametric form
\begin{equation}\label{eq:IntOpParametric}
  (\opI^q \varphi)(x) = \int_{[-1,1]^2} K^q(x,u,v) \varphi^q(u,v) J^q(u,v) \, dudv, 
\end{equation}
for $x \in \Gamma$, where $K^q(x,u,v) = K(x,y^q(u,v))$ and
$\varphi^q(u,v) = \varphi(y^q(u,v))$, and where $J^q(u,v)dudv$ denotes
the element of area. We refer to $x$ as an ``evaluation'' or
``target'' point. In what follows, the
``interaction''~\eqref{eq:IntOpParametric} between an integration
patch $\Gamma^q$ and a target point $x$ will be said to be singular,
near-singular, or regular, depending on whether the distance between
the point $x$ and the integration patch $\Gamma^q$ is less than or
greater than a certain ``proximity distance'' $\delta > 0$.  Let
\begin{equation}\label{eq:PointSetDist}
  \dist{ x, \Gamma^q } = 
    \inf \left\{\, |x - y| \,\, | \,\, y \in \Gamma^q \right\}
\end{equation}
denote the distance from a point $x$ to a patch $\Gamma^q$ (where
$|\cdot|$ denotes the Euclidean distance). Then, the set of
``singular'' and ``nearly-singular'' target points associated with
$\Gamma^q$ is defined by
\begin{equation}\label{eq:SingGridPts}
  \Omega_q^{s,\delta} = \left\{ x \in \Gamma \,\, | \,\, 
                        \dist{ x, \Gamma^q } \leq \delta \right\}.
\end{equation}
In contrast, the set of regular (non-singular) target points is
defined by
\begin{equation}\label{eq:ReguGridPts}
  \Omega_q^{r,\delta} = \left\{ x \in \Gamma \,\, | \,\, 
                        \dist{ x, \Gamma^q } > \delta \right\}.
\end{equation}
Interactions between $\Gamma^q$ and $x\in\Gamma$ with
$x\in \Omega_q^{s,\delta}$ are said to be
\emph{singular/near-singular}; other interactions, wherein
$x\in\Omega_q^{r,\delta}$,
are said to be \emph{regular}, or,
alternatively, \emph{non-singular}.

The accurate evaluation of the singular and near singular integrals
require special treatment, as discussed in
Section~\ref{sec:SingInter}, on account of the correspondingly
singular or near-singular character of the operator kernels in such
cases. The regular integrations, in turn, do not require use of
sophisticated quadrature methods, but they involve integration over
most of the surface $\Gamma$---and thus give rise to a very large
computational costs, for high-frequency problems, if straightforward
quadrature methods are used. To avoid the latter problem, the IFGF
method introduced in the previous section is utilized in this paper:
as shown in Section~\ref{sec:SurfRep}, all regular integrations may be
expressed in the form~\eqref{sec:IFGF} for whose evaluation the IFGF
method was designed.

The quadrature methods used in this paper are based on representation
of $\varphi^q$ via Chebyshev expansion for each $q=1,\dotsc,Q$. To
incorporate the discrete Chebyshev framework in our context we
discretize each patch $\Gamma^q$ by means of a surface grid
$\Bdy^q_{N^q_u,N^q_v}$ containing $N_u^q \times N_v^q$ discretization
points obtained as the image of the tensor-product discretization
\[
  \left\{ u_i = s^q_i \,\,|\,\, i = 0,\dotsc, N_u^q-1 \right\} \times
  \left\{ v_j = t^q_j \,\,|\,\, j = 0,\dotsc, N_v^q-1 \right\},
\]
under the surface parametrization $y^q$ in equation~\eqref{param_def}, where
$s^q_i$ and $t^q_j$ denote the the nodes Chebyshev-Gauss nodes
\begin{align}\label{eq:FejersNodes}
  s^q_i = \cos \left( \pi \frac{2i+1}{2 N^q_u} \right), \qquad
  t^q_j = \cos \left( \pi \frac{2j+1}{2 N^q_v} \right),
\end{align}
for $i = 0,\dotsc,N^q_u-1$, $j = 0,\dotsc,N^q_v-1$, and $q = 1,\dotsc,Q$.  The
set of all surface discretization points is denoted by
\begin{equation}\label{eq:AllSurfGrid}
  \Bdy_N = \bigcup_{q=1}^Q \Bdy^q_{N^q_u,N^q_v},
\end{equation}
where $N$ denotes the total number of grid points over all component patches.

For $1\leq q\leq Q$, the surface density $\varphi^q$ admits the
Chebyshev approximation
\begin{equation}\label{eq:ChebyExpansion}
  \varphi^q(u,v) \approx \sum_{m=0}^{N_v^q-1} \sum_{n=0}^{N_u^q-1}
                         a_{n,m}^q T_n(u) T_m(v),
\end{equation}
where, in view of the discrete orthogonality property satisfied by
Chebyshev polynomials, we have
\begin{equation}\label{eq:ChebyCoeff}
  a_{n,m}^q = \frac{\alpha_n \alpha_m}{N_u^q N_v^q} 
  \sum_{j=0}^{N_v^q-1} \sum_{i=0}^{N_u^q-1} \varphi^q(s^q_i,v_j) 
  T_n(s^q_i) T_m(v_j), 
  \qquad \alpha_n = \begin{cases} 1, & n=0 \\ 
                              2, & n \neq 0.
                                        \end{cases}
\end{equation}

\subsection{Non-singular integration via IFGF
  acceleration}\label{sec:ReguInter}
To evaluate the integral operator~\eqref{eq:IntOpParametric} at
regular target points $x_{\ell} \in \Omega_q^{r,\delta}$, we use
Fej\'{e}r's first quadrature rule, which is based on the
nodes~\eqref{eq:FejersNodes} and the integration weights
\begin{align}\label{eq:FejersWeights}
 w_j =  w^J_j &= \frac{2}{J} \left[1 - 2\sum_{\ell = 1}^{\floor{J/2}} \frac{1}{4\ell^2 -1} 
         \cos \left( \ell \pi \frac{2j+1}{J} \right) \right], 
\end{align}
so that the value at $x_{\ell}$ is approximated by
\begin{equation}\label{eq:ReguIntOpQuadrature}
  (\opI^q \varphi)(x_{\ell}) \approx 
    \sum_{j=0}^{N_v^q-1} \sum_{i=0}^{N_u^q-1} K^q(x_{\ell},s^q_i,t^q_j) 
      \varphi^q(s^q_i,t^q_j) J^q(s^q_i,t^q_j) \, w^{N_u^q}_i w^{N_v^q}_j.
\end{equation}

Evaluating all regular (non-singular) interactions on the basis
of~\eqref{eq:ReguIntOpQuadrature} results in an algorithm with
computing cost of the order of $\mathcal{O}(N^2)$ operations, whereas
the cost of adding all of the singular interactions (by means of the
RP algorithm presented in Section~\ref{sec:SingInter}) requires only
$\mathcal{O}(N)$ cost---since refinements of discretizations are
obtained by partitioning patches and keeping constant the number of
points per patch.  For acoustically-large problems such a
computational cost becomes prohibitively high.  Noting that the sum of
the regular interactions~\eqref{eq:ReguIntOpQuadrature} over all $q$
and for each $\ell$ can manifestly be expressed in the
form~\eqref{eq:field1}, to address this difficulty we accelerate the
evaluation of regular interactions by resorting to the IFGF
acceleration method introduced in Section~\ref{sec:IFGF}.  This
results in an overall iterative solver with asymptotic
$\mathcal{O}(N\log N)$ computing cost. To complete our description of
the proposed IFGF-accelerated scattering solver algorithm it therefore
suffices to introduce an accurate and efficient algorithm for singular
and near-singular interactions---which, as indicated above, is
presented in Section~\ref{sec:SingInter}.


\subsection{Singular and near-singular integration via the
  Rectangular-Polar method}\label{sec:SingInter}
To evaluate~\eqref{eq:IntOpParametric} at a singular or near-singular
target point $x\in \Omega_p^{s,\delta}$ ($\delta \geq 0$) the RP
method proceeds as follows.  Replacing the density $\varphi^q$ by its
Chebyshev expansion~\eqref{eq:ChebyExpansion}, in the RP method the
integral~\eqref{eq:IntOpParametric} is numerically approximated by
\begin{subequations}\label{eq:IOCheExpansion}
\begin{align}
  (\opI^q \varphi)(x) &\approx \int_{R} K^q(x,u,v)   \left( \sum_{m=0}^{N_v^q-1} \sum_{n=0}^{N_u^q-1} a_{n,m}^q T_n(u) T_m(v) \right)  J^q(u,v) \, dudv \label{eq:IOCheExpA} \\
 &= \sum_{m=0}^{N_v^q-1} \sum_{n=0}^{N_u^q-1} a_{n,m}^q   \left( \int_{R} K^q(x,u,v)   T_n(u) T_m(v) J^q(u,v) \, dudv \right). \label{eq:IOCheExpB}
\end{align}
\end{subequations}
Note that the double integral in this equation does not depend on the
density $\phi^q$: it depends only on the kernel, a product of
Chebyshev polynomials, and the geometry. The proposed method proceeds
by accurately precomputing and storing such quantities for all
relevant discretization points $x=x_\ell$; the algorithm is completed
by direct evaluation of the sum~\eqref{eq:IOCheExpB}. In other words,
letting
\begin{equation}\label{eq:BetaAtSing}
  \beta_{n,m}^{q,\ell} = \int_{R} K^q(x_{\ell},u,v) 
                                   T_n(u) T_m(v) J^q(u,v) \, dudv,
\end{equation}
the algorithm obtains the value of $\opI^q$ at all target points
$x_{\ell} \in \Omega_p^{s,\delta}$ on the basis of the expression
\begin{equation}\label{eq:IntOpAtSing}
  (\opI^q \varphi)(x_{\ell}) = \sum_{m=0}^{N_v^q-1} \sum_{n=0}^{N_u^q-1} a_{n,m}^q 
                                 \, \beta_{n,m}^{q,\ell}.
\end{equation}

To obtain~\eqref{eq:BetaAtSing} at an evaluation point
$x_{\ell}\in\Gamma$ the algorithm utilizes the point
$x^q_{\ell}\in \Gamma^q$ that is closest to $x_{\ell}$, or, more
precisely, the parameter-domain coordinates of $x^q_{\ell}$ in the
patch $\Gamma^q$.  If $x^q_{\ell}$ is itself a grid point of
$\Gamma^q$,
$x^q_{\ell} = x_{\ell} =y^q(\bar{u}_{\ell}^q, \bar{v}_{\ell}^q)$, then
its parameter-space coordinates in the patch $\Gamma^q$ are of course
equal to $(\bar{u}_{\ell}^q, \bar{v}_{\ell}^q)$.  More generally,
$(\bar{u}_{\ell}^q, \bar{v}_{\ell}^q)$ is obtained as the distance
minimizer:
\begin{equation}\label{eq:MinNode}
  (\bar{u}_{\ell}^q, \bar{v}_{\ell}^q) = \argmin_{(u,v)\in [-1,1]^2}
                                         \norm{ x_{\ell} - y^q(u,v) }.
\end{equation}
Following~\cite{BrunoGarza2020}, for robustness and simplicity our
algorithm tackles the minimization problem~\eqref{eq:MinNode} by means
of the golden section search algorithm.

Next, we apply a one-dimensional change of variables to each
coordinate in the $uv$-parameter space to construct a clustered grid
around each given target node.  To this end we consider the following
one-to-one, strictly monotonically increasing, and infinitely
differentiable function $w:[0,2\pi] \to [0,2\pi]$, with parameter
$d \geq 2$ (proposed in~\cite[Section 3.5]{ColtonKress2013InvAcouEM}):
\begin{equation}\label{eq:wMap}
  w( \tau; d ) = 2\pi \frac{[\nu(\tau)]^d}
                                   {[\nu(\tau)]^d + [\nu(2\pi - \tau)]^d},  
                                   \qquad 0 \leq \tau \leq 2\pi,
\end{equation}
where
\begin{equation}\label{eq:vMap}
  \nu( \tau; d ) = \left( \frac{1}{d} - \frac{1}{2} \right)
                           \left( \frac{\pi - \tau}{\pi} \right)^3
                         + \frac{1}{d} \left( \frac{\tau - \pi}{\pi} \right)
                         + \frac{1}{2}.
\end{equation}
It can be shown that $w$ has vanishing derivatives up to order $d-1$ at the
interval endpoints. Then, the following change of variables 
\begin{equation}\label{eq:CoV}
  \xi_{\alpha}( \tau; d ) = 
  \begin{cases}
  \alpha + \left( \frac{\sgn(\tau) - \alpha}{\pi} \right) w( \pi |\tau|; d ), 
    & \text{ for } \alpha \neq \pm 1, \\
  \alpha - \left( \frac{1+\alpha}{\pi} \right) w( \pi \left| \frac{\tau - 1}{2} \right|; d ), 
    & \text{ for } \alpha =  1, \\
  \alpha + \left( \frac{1-\alpha}{\pi} \right) w( \pi \left| \frac{\tau + 1}{2} \right|; d ), 
    & \text{ for } \alpha = -1,
  \end{cases}
\end{equation}
has the effect of clustering points around $\alpha$.  Fej\'{e}r's rule
(Section~\ref{sec:ReguInter}) applied to the
integral~\eqref{eq:BetaAtSing}, transformed using the change of
variables~\eqref{eq:CoV}, yields the approximation
\begin{equation}\label{eq:BetaCov}
  \beta_{n,m}^{q,\ell} \approx \sum_{j=0}^{N_{\beta}-1} \sum_{i=0}^{N_{\beta}-1}
                               K^q(x_{\ell}, u_i^{q,\ell}, v_j^{q,\ell} ) 
                               T_n(u_i^{q,\ell}) T_m(v_j^{q,\ell}) 
                               J^q(u_i^{q,\ell},v_j^{q,\ell}) \,
                               w_i^{u,q,\ell} w_j^{v,q,\ell},
\end{equation}
where
\begin{align}\label{eq:CovNodesWeights}
  u_i^{q,\ell} &= \xi_{\bar{u}_{\ell}^q} ( s_i; d ), \quad
  w_i^{u,q,\ell} = \frac{ d\xi_{\bar{u}_{\ell}^q} }{d\tau}( s_i; d ) \, w_i, \\
  v_j^{q,\ell} &= \xi_{\bar{v}_{\ell}^q} ( s_j; d ), \quad
  w_j^{v,q,\ell} = \frac{ d\xi_{\bar{v}_{\ell}^q} }{d\tau}( s_j; d ) \, w_j,
\end{align}
for $i,\, j = 0, \dotsc, N_{\beta} - 1$.  To avoid division by zero,
we set the kernel $K^q$ to zero at integration points where the
distance to the target point is less than some prescribed tolerance,
usually on the order of $10^{-14}$. This completes the singular and
near-singular integration algorithm, and, thus, incorporating the
iterative linear algebra solver GMRES, the proposed overall
IFGF-accelerated high-order integral solver for the
problem~\eqref{eq:AcousticBVP}.

\section{Numerical results}\label{sec:NumRes}

This section presents numerical results that demonstrate the character
of the proposed IFGF-accelerated acoustic scattering solvers described
in previous sections.  For comparison, results obtained using the
$\Ord(N^2)$ non-accelerated Chebyshev-based scattering solvers
introduced in~\cite{BrunoGarza2020} are also included. Both the
accelerated and non-accelerated solver are implemented using OpenMP for
shared-memory parallelism.

After solving~\eqref{eq:CombLayIE} for the density $\varphi$, the far field
pattern $u^{\infty}$ can be obtained from 
\begin{equation}\label{eq:FarField}
  u^{\infty}(\hat{x}) = \frac{1}{4\pi} 
    \int_{\Bdy} \left\{ \frac{\del}{\del \nu(y)} e^{-\ii k \hat{x} \cdot y} 
                                  - \ii \gamma e^{-\ii k \hat{x} \cdot y} \right\} 
      \varphi(y) \, dS(y), \quad \hat{x} \in \mathbb{S}^2,
\end{equation}
where $\mathbb{S}^2$ denotes the unit sphere and $\Bdy$ is the scatterer's
boundary. The far field is computed over a uniformly-spaced unit spherical grid
\begin{equation}\label{eq:SphereGrid}
  \Sfar = \left\{ (\phi_m,\theta_n) \in [0,\pi]\times[0,2\pi] 
                              \,\,|\,\, 1 \leq m \leq N_{\phi}, \, 
                                        1 \leq n \leq N_{\theta} \right\},
\end{equation}
with $\phi_m = (m-1)\Delta \phi, \, \theta_n = (n-1) \Delta \theta$
and where the spacings are defined as $\Delta \phi = \pi/(N_{\phi}-1)$
and $\Delta \theta = 2\pi/(N_{\theta}-1)$, respectively; specific
values of $N_{\phi}$ and $N_{\theta}$ are given in each example's
subsection.  Given the discrete values of exact (reference) far field
modulus $|u^{\infty}_{m,n}|$ and approximate far field modulus
$|\tilde{u}^{\infty}_{m,n}|$
($1 \leq m \leq N_{\phi}, \, 1 \leq n \leq N_{\theta}$), the maximum
far field relative error $\farErr$ over $\Sfar$ given by
\begin{equation}\label{eq:FarFieldRelErr}
  \farErr = \max_{(m,n) \in \Sfar} \left\{
    \frac{| \ |u_{m,n}^{\infty}| -  |\tilde{u}_{m,n}^{\infty}|\ |}
    { |u_{m,n}^{\infty}| } \right\},
\end{equation}
is reported in each case.

Similarly, substituting the solution $\varphi$ in the combined-layer
representation~\eqref{eq:CombinedLayer} we evaluate and display the
scattered field $u^s$ over near field planes that are parallel to the
$xy$-, $xz$-, or $yz$-planes.  For example, we evaluate fields
(incident, scattered, and total) at every point of a uniformly-spaced
two-dimensional $xy$-planar grid $\Pnear^{xy}(z_0)$ at $z = z_0$
defined by
\begin{equation}\label{eq:PlanarGrid}
  \Pnear^{xy}(z_0) = \left\{ (x_m,y_m,z_0) \in [x_{min},x_{max}] \times 
    [y_{min},y_{max}]  \times \{ z_0 \}
                              \,\,|\,\, 1 \leq m \leq N_x, \, 
                                        1 \leq n \leq N_y \right\},
\end{equation}
where the grid points are given by
$x_m = (m-1) \Delta x, \, y_n = (n-1) \Delta y$ and the grid spacings
are $\Delta x = (x_{max} - x_{min}) /(N_x-1)$ and
$\Delta y = (y_{max} - y_{min}) /(N_y-1)$. Near field planar grids
parallel to the $xz$- and $yz$-plane are defined analogously. Denoting
the exact (or reference) and approximate modulus of the total field at
each point of $\Pnear^{xy}(z_0)$ by $v_{m,n}$
($=|u_{m,n}^s+u^i_{m,n}|$) and $\tilde{v}_{m,n}$
($=|\tilde{u}_{m,n}^s+u^i_{m,n}|$), respectively, we define the near
field (total magnitude) relative error $\nearErr$ over
$\Pnear^{xy}(z_0)$ by
\begin{equation}\label{eq:NearFieldRelErr}
  \nearErr = \max_{(m,n) \in \Pnear^{xy}(z_0)} 
    \left\{ \frac{|v_{m,n} - \tilde{v}_{m,n}|} {|v_{m,n}|} \right\}.
\end{equation}

The numerical results presented in what follows were obtained using 28
cores in a single Intel Xeon Platinum 8276 2.20 GHz computer. The
images presented were generated using the visualization software
VisIt~\cite{HPV:VisIt}. Solutions to the complex-coefficient linear
systems that arise from discretizations of the boundary integral
equation~\eqref{eq:CombLayIE} were obtained by means of a
complex-arithmetic GMRES iterative solver~\cite{Saad1986}.
Following~\cite{BrunoKunyansky2001JCP}, we utilize the value
\begin{equation}
  \label{eq:gamma}
  \gamma = \max\{ 3, A/\lambda \}
\end{equation}
of the coupling parameter in the combined-layer
equation~\eqref{eq:CombLayIE}, where $A$ denotes the diameter of the
scatterer. Detailed computational studies indicate that, to reach a
given residual tolerance, this value reduces the number of GMRES
iterations by a factor of $5-10$ vs. the iteration numbers required
when the oft considered value $\gamma = k$ is used. Per
Remark~\ref{cone_seg_size}, in all cases the values $n_{C,0}=2$ and
$n_{S,0} = 1$ were utilized.

\subsection{Scattering by a sphere}\label{sph_scat}

%
%
This section concerns the classical problem of scattering of plane
waves by sound-soft spheres, for which the well-known closed-form Mie
series far field expression is used to compute relative
errors~\cite{BowmanEtAl1988}. Spheres of various acoustical sizes are
considered, ranging from $4\lambda$ to $128\lambda$ in diameter. In
all cases the number of IFGF levels was selected so that the
finest-level IFGF box side length is approximately $0.5\lambda$.

\begin{figure}[h!]
    \centering
    \includegraphics[width=0.7\textwidth]{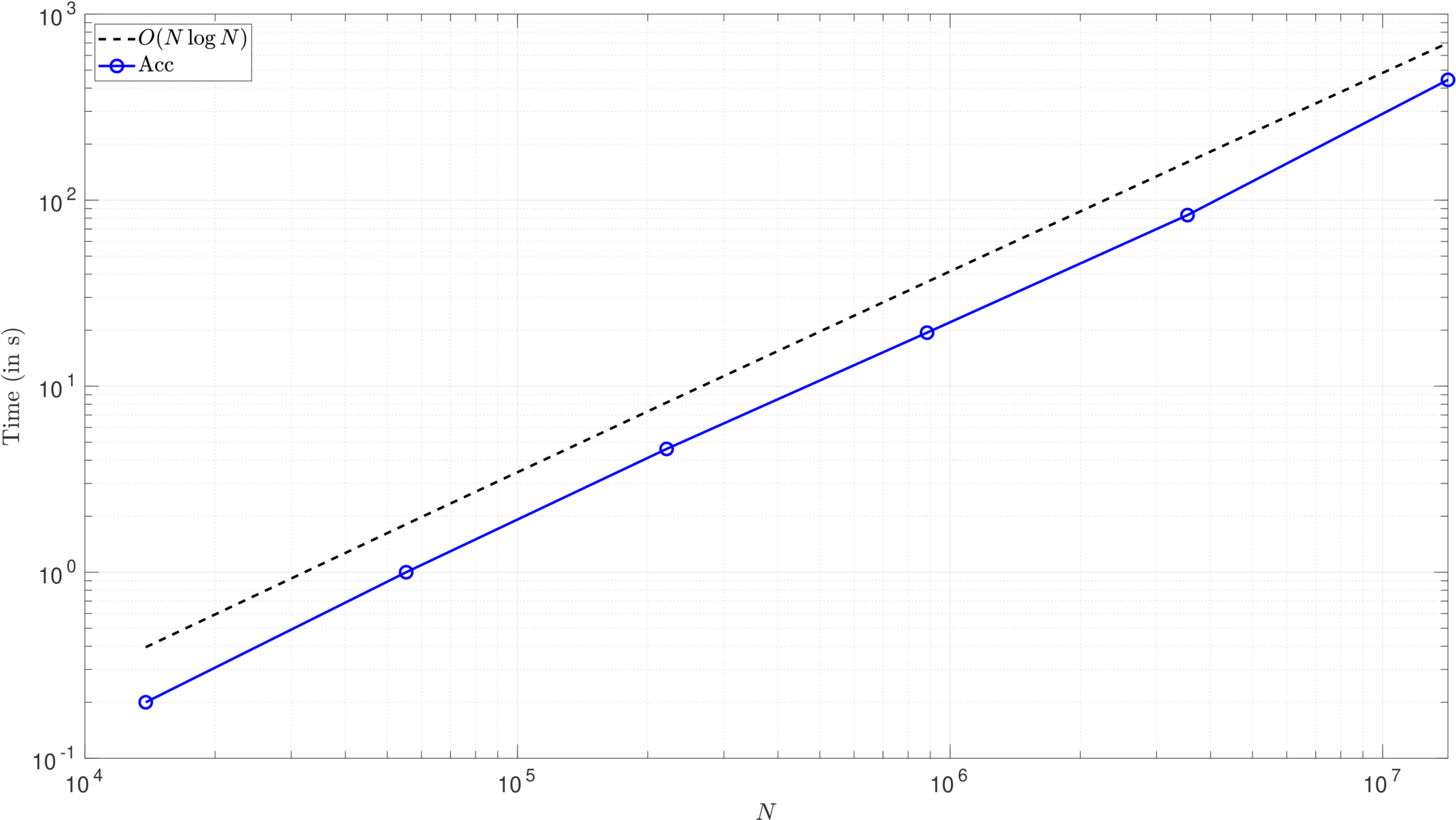}
    \caption{\small Computing times, in seconds, required by the
      IFGF-accelerated solver, as a function of the total number $N$
      of discretization points (line-circle markers), for the problem
      of acoustic scattering by spheres.  For reference, the figure
      includes a curve parallel the graph of $C N \log N$ with slope
      $C=0.3 \cdot 10^{-5}$ (dashed line).}
  \label{fig:SphereNLogN}
\end{figure}

Table~\ref{table:AcouSphere} demonstrates the character of both, the
IFGF-accelerated solver and the non-accelerated solver, via
applications to problems of scattering by spheres of diameters ranging
from $4$ to $128$ wavelengths.  In all cases the GMRES residual
tolerance was set to $10^{-4}$. We report the total number of
unknowns, the size of the sphere in wavelengths, the time required to
compute one GMRES iteration, the total number of iterations
required to achieve the prescribed residual, and the field relative
error. Far field relative errors are presented over the spherical
grid~\eqref{eq:SphereGrid} with $(N_{\phi},N_{\theta}) =
(200,200)$. The computing times per iteration displayed in this table
for the non-accelerated algorithm, grow by a factor of around
$14.8-15.7$ as the number of points per dimension in each surface
patch is doubled (so that the overall number $N$ of unknowns is
quadrupled), which is consistent with the expected quadratic
complexity of the algorithm. The corresponding costs for the
IFGF-accelerated solver, on the other hand, scale like
$\Ord(N\log N)$, as illustrated in Figure~\ref{fig:SphereNLogN}. In
all cases, solutions with errors of the order of $10^{-4}$ were
obtained.

\begin{table}[H] 
\centering
\begin{tabular}{lcccccccc}\toprule
  & & \multicolumn{3}{c}{Nonaccelerated} & \multicolumn{4}{c}{Accelerated} 
  \\ \cmidrule(lr){3-5} \cmidrule(lr){6-9}
  Unknowns    &Size         &Time (1 iter.)  &Iters. &$\farErr$ 
                &Levels  &Time (1 iter.)  & Iters. &$\farErr$  \\ 
  \midrule
  13,824      &4$\lambda$   &0.5 s           &12         & $1.1\cdot 10^{-4}$   
                &4            &0.2 s           &12         & $1.3\cdot 10^{-4}$ \\
  55,296      &8$\lambda$   &7.4 s           &14         & $8.9\cdot 10^{-5}$ 
                &5            &1.0 s           &14         & $1.1\cdot 10^{-4}$ \\
  221,184     &16$\lambda$  &116.4 s         &14         & $2.6\cdot 10^{-5}$ 
                &6            &4.6 s           &14         & $6.3\cdot 10^{-5}$ \\
  884,736     &32$\lambda$  &1862.4 s (est.) &$\--$      &$\--$  
                &7            &19.4 s          &16         & $2.9\cdot 10^{-5}$\\
  3,538,994   &64$\lambda$  &8.3 h (est.)    &$\--$      &$\--$    
                &8            &83.1 s          &18         & $6.0\cdot 10^{-5}$ \\ 
  14,155,776  &128$\lambda$  &132.8 h (est.) &$\--$      &$\--$    
                &9            &443.2 s          &21         & $3.8\cdot 10^{-4}$ \\ 
 \bottomrule
\end{tabular}
\caption{\small Performance of both, the IFGF-accelerated solver and the
  non-accelerated solver, for problems of scattering by spheres of
  diameters ranging from $4$ to $128$ wavelengths. The table
  summarizes the total number of surface unknowns, sphere size in
  wavelengths, number of IFGF levels used, time required to compute
  one GMRES iteration, total number of iterations, and far field
  relative error $\farErr$.  In all cases the GMRES residual tolerance
  was set to $10^{-4}$.}
\label{table:AcouSphere}
\end{table}
The reduced complexity of the IFGF-based algorithm has a significant
impact on computing times. At $128$ wavelengths, the computing time
per iteration required by the non-accelerated solver are over $1000$
times larger than those required by the accelerated method; for larger
problems, the difference in compute times grows as expected from the
complexity estimates for the two methods. It is worthwhile to note
that, in all cases for which iteration numbers are presented for both
methods, the non-accelerated and accelerated solvers require the same
number of GMRES iterations to meet the prescribed residual tolerance.
Additionally, the errors for both algorithms are comparable: the
non-accelerated solver yields solutions for the $4,\,8$ and $16$
wavelength problems with an average relative error of
$1.1 \cdot 10^{-4}$, while errors obtained with the accelerated method
average $1.3 \cdot 10^{-4}$ across the entire $4$ to $128$ wavelength
range.

\subsection{Comparisons with FMM-based solvers}\label{compari}

This section provides comparisons of the proposed IFGF/RP algorithm
with two recent FMM accelerated scattering solvers for the problem
considered in this paper. The first test case considered, whose
results are presented in Table~\ref{tab:ErrorVsGMRES}, examines the
convergence of the far field scattered by a spherical scatterer as a
function of the number of iterations in the GMRES iterative solver,
for the formulation~\eqref{eq:CombLayIE} and with the value of
$\gamma$ given in~\eqref{eq:gamma}. The selected problem, a sphere
$2$m in diameter under illumination at the frequency $f = 14.880$ KHz,
coincides with the third test case considered
in~\cite[Table~2]{AbduljabbarEtAl2019} for the single-layer,
FMM-accelerated formulation introduced in that reference. The present
results compare favorably with the one provided in that reference:
using 1,572,864 degrees of freedom (DoF), a total of 8 iterations of
the present solver suffice to yield a far-field accuracy of
$1.3\cdot 10^{-3}$ whereas 1,123 iterations and $8,984,640$ DoFs are
utilized in that reference to obtain a far-field error of
$3\cdot 10^{-3}$. Clearly, the second-kind
formulation~\eqref{eq:CombLayIE} utilized in this paper is greatly
advantageous, in spite of the fact of its requirement of evaluation of
two operators (the single- and double-layer potentials) instead of one
operator in the single-layer formulation used
in~\cite{AbduljabbarEtAl2019}. The computing times per iteration
required by IFGF/RP algorithm to achieve this accuracy also compare
favorably: each two-operator (single and double layer) iteration of
the present solver requires a computing time of $58.3$s in a total of
28 cores, whereas the one-operator (single layer) FMM-based iteration
in~\cite{AbduljabbarEtAl2019} requires $46.8$s in a total of 512
cores.
\begin{table}[h!]
    \centering
    \begin{tabular}{|c|cccc|}
      \hline
      \# Iters. &  2 & 4 & 8 & 12\\
      \hline
      $\farErr$  & $1.9\cdot 10^{-1}$ & $1.8\cdot 10^{-2}$ & 1.3 $\cdot 10^{-3}$ & $5.8\cdot 10^{-4}$\\
      Run time (in s, at 58.3 s/iter) & 116.6 & 233.2 & 466.4 & 699.6\\     \hline
 \end{tabular}
 \caption{\small Far-field relative error resulting from the proposed
   accelerated solver as a function of the number of GMRES iterations
   incurred, for a configuration considered
   in~\cite[Table~2]{AbduljabbarEtAl2019}: scattering by a sphere $2$
   m in diameter at a frequency of $f = 14.88$ KHz. At the speed of
   sound assumed in that paper ($c = 343$ m/s), the sphere is
   approximately $86.8\lambda$ in diameter
   ($2\mathrm{m}/\lambda = 2\mathrm{m} f/c = 86.8$). The results
   presented in the present table were obtained on the basis of a
   total of $1,572,864$ unknowns.}
  \label{tab:ErrorVsGMRES}
\end{table}

A second point of reference is provided in the
publication~\cite{WalaKlockner2019}, which demonstrates an
FMM-accelerated QBX algorithm for a problem of scattering by an
aircraft geometry with a maximum dimension of approximately
$64\lambda$, on a 20-core 2.30 GHz Intel Xeon E5-2650 v3 machine. A
sequence of geometry-specific optimizations listed
in~\cite[Fig. 11]{WalaKlockner2019} reduces a certain reference FMM
computing times by certain percentages. The optimized, 14 million DoF
implementation is reported to run at computing times for the
double-layer and single-layer operators of 1600 seconds and 900
seconds, respectively, for a total of 2,500 seconds in the combined
formulation used, to achieve an estimated accuracy of
$2.5\cdot 10^{-3}$. Results of a run containing a similar number of
degrees of freedom on a sphere, but for the much larger acoustical
size of $128\lambda$ in diameter, are presented in the last row of
Table~\ref{tab:ErrorVsGMRES}---showing a 28-core computing time (2.20
GHz cores) of 443.2s per iteration, together with an error of
$3.8\cdot 10^{-4}$.  In all, these experiments and comparisons suggest
that the proposed IFGF/RP method provides a valuable new tool for the
solution of high frequency scattering problems.

\subsection{Scattering by a submarine geometry}\label{sub_scat}
\begin{figure}[h!]
    \centering
    \includegraphics[width=\textwidth]{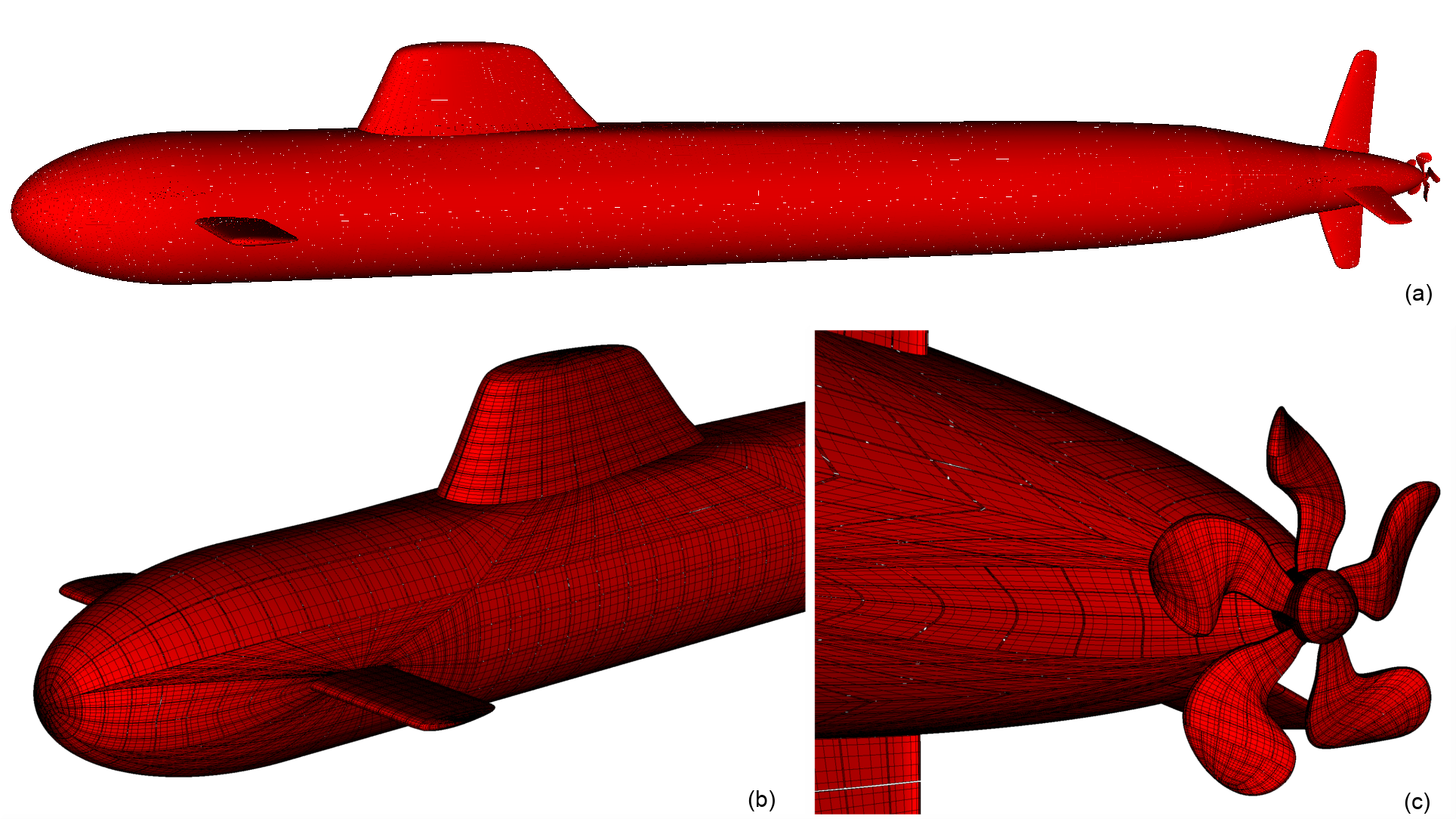}
    \caption{\small Submarine model and surface mesh containing a
      total of 164,160 points.  The submarine hull is aligned with
      the $z$-axis and the sail is parallel to the $+y$-axis; the
      front of the vessel points in the $+z$-direction.}
  \label{fig:SubmarineModel}
\end{figure}

Owing to its importance in detection and tracking applications, the
accurate simulation of problems of acoustic scattering involving
submarine vehicles is the subject of ongoing research
~\cite{NellGilroy2003,SchneiderEtAl2003,MerzEtAl2009,Karasalo2012,WeiEtAl2016,TestaGreco2018,VenaasKvamsdal2020}. This
section demonstrates the character of the proposed IFGF accelerated
solver in such a context: it concerns problems of acoustic scattering
by realistic submarine configurations of up to $80$ wavelengths in
acoustical size. The submarine model used, comprising the main hull,
sail, diving planes, rudders, and a five-blade propeller, is depicted
in Figure~\ref{fig:SubmarineModel}.  The complete submarine geometry
is contained in the bounding region
$[-3.2,3.2] \times [-1.9,2.8] \times [-19.2,10.9]$.
Figures~\ref{fig:SubmarineModel}(b) and~\ref{fig:SubmarineModel}(c)
display portions of a surface representation containing $4,560$
patches, each of which is discretized by $6\times 6$ points.

We consider plane wave scattering for two cases: a) head-on incidence and b)
oblique incidence.  The incident field is a plane wave $u^i$ that travels along
the wave direction $\hat{k}$ and is given by
\begin{equation}\label{eq:PlaneWave}
  u^i(x) = e^{ \ii k \hat{k} \cdot x}, \qquad 
  \hat{k} = \begin{pmatrix} 
              \cos \theta \sin \phi \\ 
              \sin \theta \sin \phi \\  
              \cos \phi 
            \end{pmatrix},
\end{equation}
where $x = (x_1,x_2,x_3)$ denotes the position vector, $k>0$ is a
given wavenumber, and $(\theta,\phi) \in [0,2\pi) \times
[0,\pi]$. Since the bow of the submarine points in the $+z$-direction,
``head-on'' incidence corresponds to $(\theta,\phi) = (0,\pi)$
in~\eqref{eq:PlaneWave}. For the oblique incidence case we set
$(\theta,\phi) = (0,5\pi/4)$.
%
%

To verify the accuracy of the IFGF-accelerated solver in the present
case, we conducted convergence studies for the submarine structures of
acoustical sizes (measured from the bow to the propeller cap) of
$10\lambda,\,20\lambda$, $40\lambda$ and $80\lambda$. In all cases the
number of IFGF levels was chosen so that the side length of the
smallest, (finest level) boxes is approximately $0.8\lambda$, and the
GMRES residual tolerance was set to $10^{-3}$. We start our sequence
of test cases with a $10\lambda$ vessel whose geometry is represented
by $1,140$ surface patches, each of which has $6 \times 6$ points.  As
the size of the problem is doubled, the geometry is partitioned from
the previous size so that every patch is split into four subpatches
while keeping the number of points per patch fixed. Thus, for example,
the $20\lambda$ problem uses four times as many surface points as the
$10\lambda$ one. This is admittedly a suboptimal strategy in this case
(since, after a certain level $d$, the smaller patches on the
propeller, rudders and diving planes are the only ones that would
provide acceleration via further partitioning), which, however,
simplifies the code implementation (cf. the discussion concerning
adaptivity at the end of Section~\ref{sec:CombLayIFGF}).  As indicated
in~\cite{BauingerBruno2021}, the IFGF method can be extended to
incorporate a box octree algorithm that adaptively partitions a
geometry so that boxes containing a (small) prescribed number of
points are not further partitioned, which eliminates the difficulty
posed in the problem at hand by finely discretized structures such as
the propeller, rudders and diving planes. While such an addition is
left for future work, as demonstrated in Table~\ref{table:Submarine},
even the simple uniform-partition IFGF algorithm we use here is
sufficient to simulate scattering by a realistic submarine geometry
for up to $80$ wavelengths in size with several digits of accuracy and
using only modest computational resources. For example, the $656,640$
unknowns, $40\lambda$ run for head-on incidence, required a computing
time of $313$ seconds per iteration and a total of $78$
iterations. The fully adaptive version of the IFGF algorithm, which,
as mentioned above, is not pursued in this paper, should yield
submarine-geometry computing times consistent with those shown in
Tables~\ref{table:AcouSphere} and~\ref{table:Nacelle} for the sphere
and nacelle geometries, respectively.
\begin{table}[htb!] 
\centering
\begin{tabular}{lcccc}\toprule
  & & & \multicolumn{1}{c}{Front Incidence} & \multicolumn{1}{c}{Oblique Incidence} 
  \\ \cmidrule(lr){4-4} \cmidrule(lr){5-5}
  Unknowns   &Size         &IFGF levels  &$\nearErr$ &$\nearErr$  \\ 
  \midrule
  41,040     &10$\lambda$  &5            & $2.4\cdot 10^{-4}$  & $6.8\cdot 10^{-4}$ \\
  164,160    &20$\lambda$  &6            & $1.8\cdot 10^{-4}$   & $6.0\cdot 10^{-4}$ \\
  656,640    &40$\lambda$  &7            & $2.1\cdot 10^{-4}$   & $1.7\cdot 10^{-4}$ \\
  2,626,560  &80$\lambda$  &8            & $2.1\cdot 10^{-4}$ (est.) & $4.8\cdot 10^{-4}$ (est.) \\
 \bottomrule
\end{tabular}
\caption{\small Near field errors resulting from applications of the
  IFGF accelerated acoustic solver to the submarine geometry depicted
  in Figure~\ref{fig:SubmarineModel}, with acoustical sizes ranging
  from $10$ to $80$ wavelengths.  Errors were evaluated by means of
  convergence studies, as detailed in the text, except for the
  $80$-wavelength case, for which, as mentioned in the text, an
  estimate was used. In all cases the residual tolerance for the
  integral density was set to $10^{-3}$.}
\label{table:Submarine}
\end{table}

Near field relative errors~\eqref{eq:NearFieldRelErr} for front
(head-on) and oblique plane wave incidence are presented in
Table~\ref{table:Submarine}. For each problem, the near field errors
$\nearErr$ were obtained over the planar region $\Pnear^{xy}(z_0)$
(equation~\eqref{eq:PlanarGrid}) with
$[x_{min},x_{max}] \times [y_{min},y_{max}] = [-12,12]^2,\, z_0 = -25$
and $N_x = N_y = 260$. The reference solution was obtained with the
same number of surface patches as the target discretization but using
$8 \times 8$ points per patch and a residual tolerance of
$10^{-5}$. (Thus, the reference solution uses nearly twice as many
discretization points and it satisfies a more stringent convergence
condition.) The numerical results indicate that the solution accuracy
is consistent for both front and oblique incidence and for all
acoustical sizes considered. The relative errors for the
$10\lambda,\, 20\lambda$, and $40\lambda$ for the front and oblique
incidence problems achieve average accuracies of $2.1\cdot 10^{-4}$
and $4.8 \cdot 10^{-4}$, respectively; these values where used to
estimate the expected relative errors in the $80$-wavelength case by
simple averaging.
\begin{figure}[h!]
  \centering
  \includegraphics[width=\textwidth]{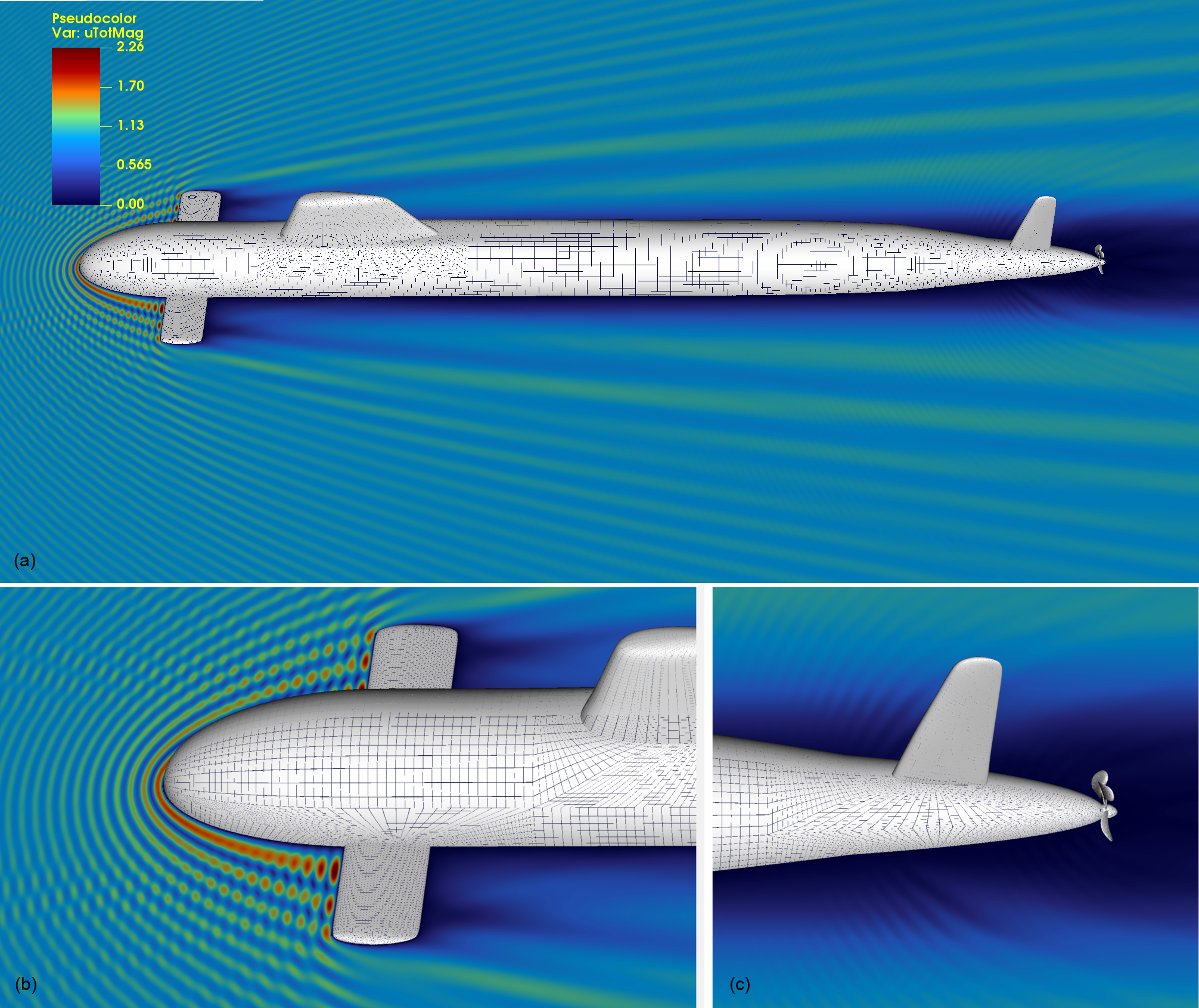}
  \caption{\small Total field magnitude $|u(x)| = |u^i(x) + u^s(x)|$
    pseudo-color plots for an $80$-wavelength submarine. The field is
    displayed over a uniform grid of $1040\times 1760$ points
    $(x,z) \in [-12,12]\times[-25,15]$.  The incident plane wave
    impinges on the vessel head-on, which corresponds to the wave
    direction $\hat{k}$ in~\eqref{eq:PlaneWave} with
    $(\theta,\phi) = (0,\pi)$.} \label{fig:SubNF0Ang}
\end{figure}
\begin{figure}[h!]
  \centering
  \includegraphics[width=\textwidth]{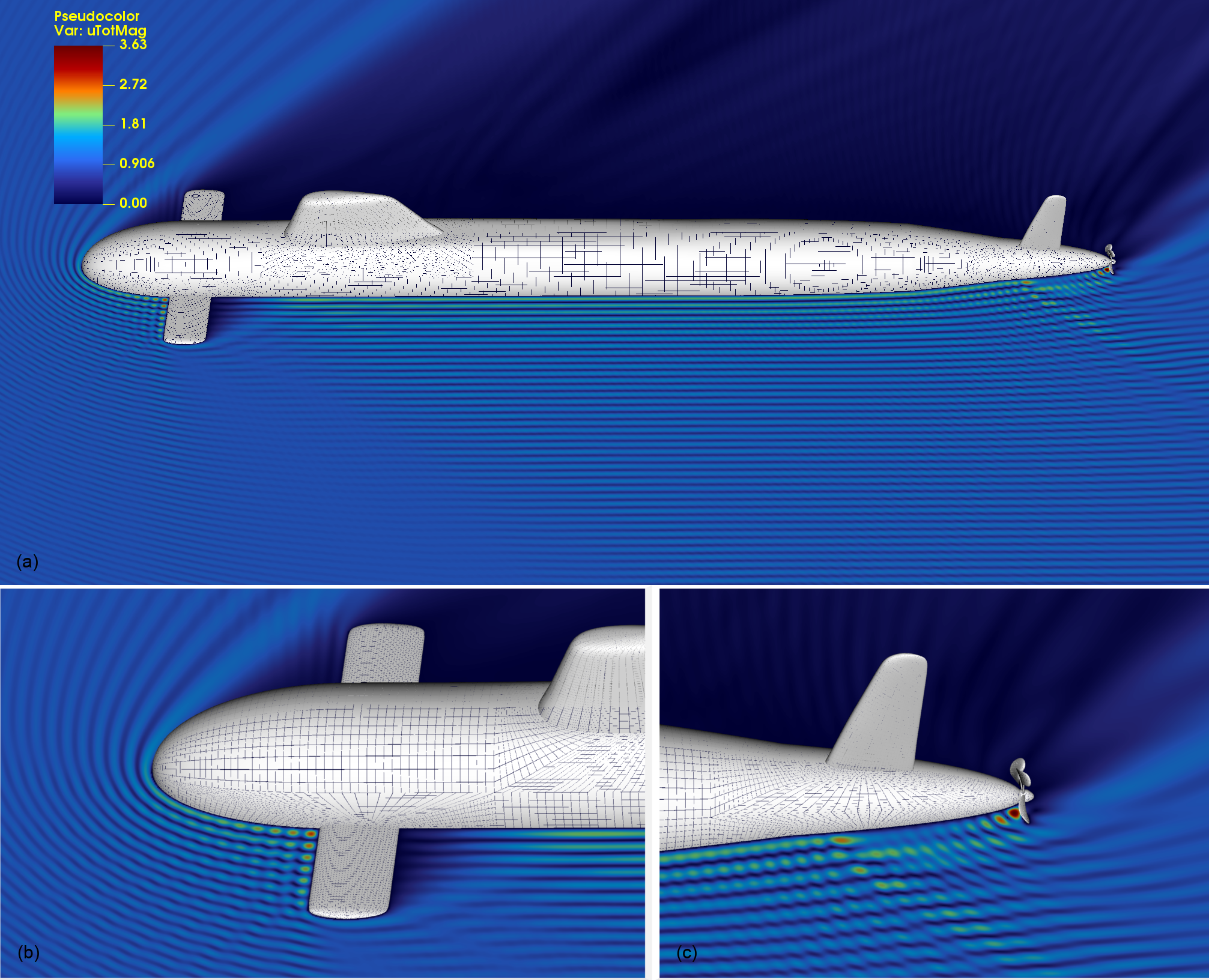}
  \caption{\small Total field magnitude $|u(x)| = |u^i(x) + u^s(x)|$
    pseudo-color plots for an $80$-wavelength submarine. The field is
    displayed over a uniform grid of $1040\times 1760$ points
    $(x,z) \in [-12,12]\times[-25,15]$.  The incident plane wave
    impinges on the vessel at an oblique angle, which corresponds to
    the wave direction $\hat{k}$ in~\eqref{eq:PlaneWave} with
    $(\theta,\phi) = (0,5\pi/4)$.}
  \label{fig:SubNF45Ang}
\end{figure}

Figure~\ref{fig:SubNF0Ang} presents pseudo-color near field plots of
the total field magnitude $|u(x)| = |u^i(x) + u^s(x)|$ resulting for
front plane wave incidence on an $80$-wavelength submarine. The field
is plotted over a uniform $1040\times 1760$ point planar grid
$\Pnear^{xz}(y^0)$ for $(x,z) \in [-12,12]\times[-25,15]$ and
$y^0 = 0$.  The incident plane wave impinges on the vessel head-on and
we see in Figures~\ref{fig:SubNF0Ang}(a) and~\ref{fig:SubNF0Ang}(b)
that the strongest interaction occurs around the bow and diving planes
(also known as hydroplanes) of the ship. Shadow regions are visible
immediately behind the hydroplanes as well as along the hull,
particularly in the aft of the ship where the body tapers.
Figure~\ref{fig:SubNF0Ang}(c) shows that the wider sections of the
ship obstruct the propeller from most incoming waves and, as a
consequence, there is minimal interaction in this region.
Figure~\ref{fig:SubNF45Ang}, in turn, shows near field pseudo-color
plots for the same $80$-wavelength submarine but this time for oblique
plane wave incidence. The total field magnitude is plotted over the
uniform grid $\Pnear^{xz}(y^0)$ described in the previous
paragraph. In this case the wave interaction is markedly different.
We see the expected shadow region in the opposite side of the incoming
wave but there is now clear evidence of wave interaction between the
hull and diving planes as well as around the rudders and propeller.
In addition to multiple scattering, the close-up views in
Figures~\ref{fig:SubNF0Ang}(b),~\ref{fig:SubNF45Ang}(b)
and~\ref{fig:SubNF45Ang}(c) show the formation of bright spots near
the junction of the left hydroplane and hull and in the vicinity of
the propeller. each one of which spawns a scattered wave propagating
in forward directions.

\subsection{Scattering by an aircraft nacelle} 
\begin{figure}[h!]
    \centering
    \includegraphics[width=\textwidth]{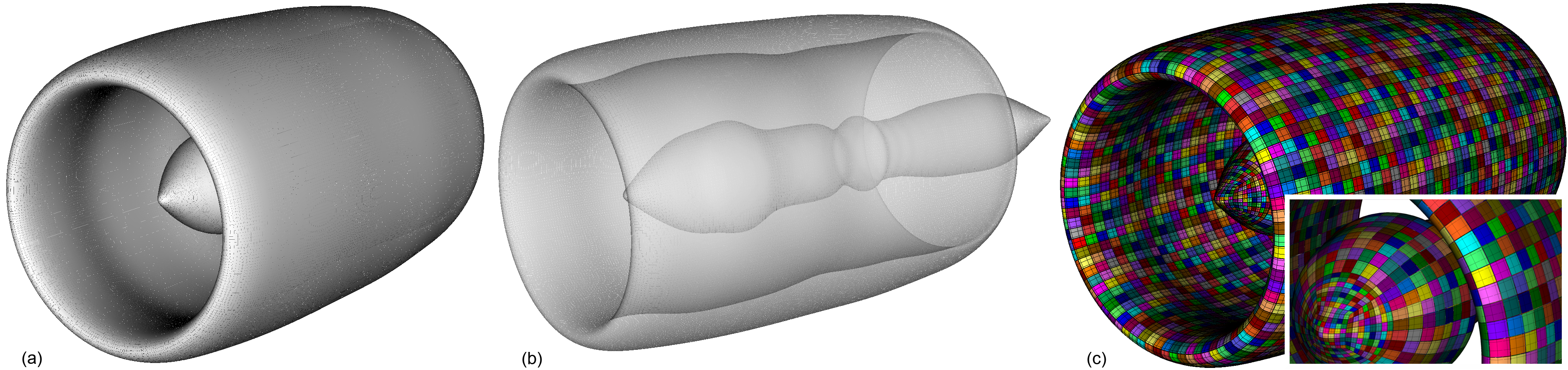}
    \caption{\small (a) Aircraft engine nacelle model, (b) translucent
      view of the geometry where the center shaft is visible, and (c)
      8,576-surface-patch discretization with $6\times6$ points per
      patch (for clarity, only every other mesh point is plotted).}
\label{fig:NacelleAndMesh}
\end{figure}

The simulation of aircraft engine noise has been the subject of
intense research for the past several decades due to its importance in
civil aviation
applications~\cite{Nayfeh1975,Eversman1991,Casalino2008,Kewin2013,PyatuninEtAl2016}.
This section presents simulations of sound scattering by the turbofan
engine nacelle model shown in Figure~\ref{fig:NacelleAndMesh} in
accordance with the boundary value
problem~\eqref{eq:AcousticBVP}. Although such acoustic simulations do
not account for many of the complexities associated with engine noise,
such as, e.g., fluid flow, density gradients, and noise generation, as
illustrated by the aforementioned references, acoustic models such as
the one considered here play significant roles in the modeling of
engine noise. According to the nacelle wall liner case
study~\cite{MSCSoft2021}, for example, under typical operating
conditions, engine nacelle noise occurs in the $125\--5650$ Hz
frequency range.  For a typical airliner engine that is around 5 m
long, these frequencies correspond to acoustical sizes between $2$ and
$82$ wavelengths. The engine nacelle geometry used in the simulations
that follow is depicted in Figures~\ref{fig:NacelleAndMesh}(a)
and~\ref{fig:NacelleAndMesh}(b); it comprises an outer housing and a
center shaft. The entire two-piece nacelle structure is contained
inside the bounding region
$[-1.5,1.5] \times [-1.5,1.5] \times [-3.27,3.27]$. The center shaft
is aligned with the $z$-axis, with the tip of the shaft pointing
towards the positive direction.  A discretization containing 8,576
surface patches with $6 \times 6$ points per patch is shown in
Figure~\ref{fig:NacelleAndMesh}(c); for future reference, note that
the inset image shows that the mesh is not rotationally symmetric,
e.g. near the tip of the shaft.

Two types of incident fields were used: a)~A plane wave that travels
towards the $-z$-axis, so that it impinges on the nacelle head-on;
and, b)~A set of eight point sources placed inside the housing around
the center shaft. As in the submarine example, the plane wave incident
field is given by~\eqref{eq:PlaneWave} with $(\theta,\phi) =
(0,\pi)$. The incident field b), on the other hand, serves as a simple
model for fan noise generation inside the nacelle and is given by
\begin{equation}\label{eq:PointSource}
  u^i(x) = \sum_{j=1}^{8} \frac{e^{\ii k |x-x^j|}}{|x-x^j|}, 
  \quad \text{with point source locations} \quad
  x^j = (x^j_1, x^j_2, x^j_3) = (\cos \alpha_j,  \sin \alpha_j, 2),
\end{equation}
where $\alpha_j = (j-1) \Delta \alpha + \pi / 8$, for
$j = 1,\dotsc, 8$, and $\Delta \alpha = \pi / 4$; as illustrated in
panel~(b) of Figure~\ref{fig:NacNFPtSrc}.
\begin{figure}[h!]
    \centering
    \includegraphics[width=\textwidth]{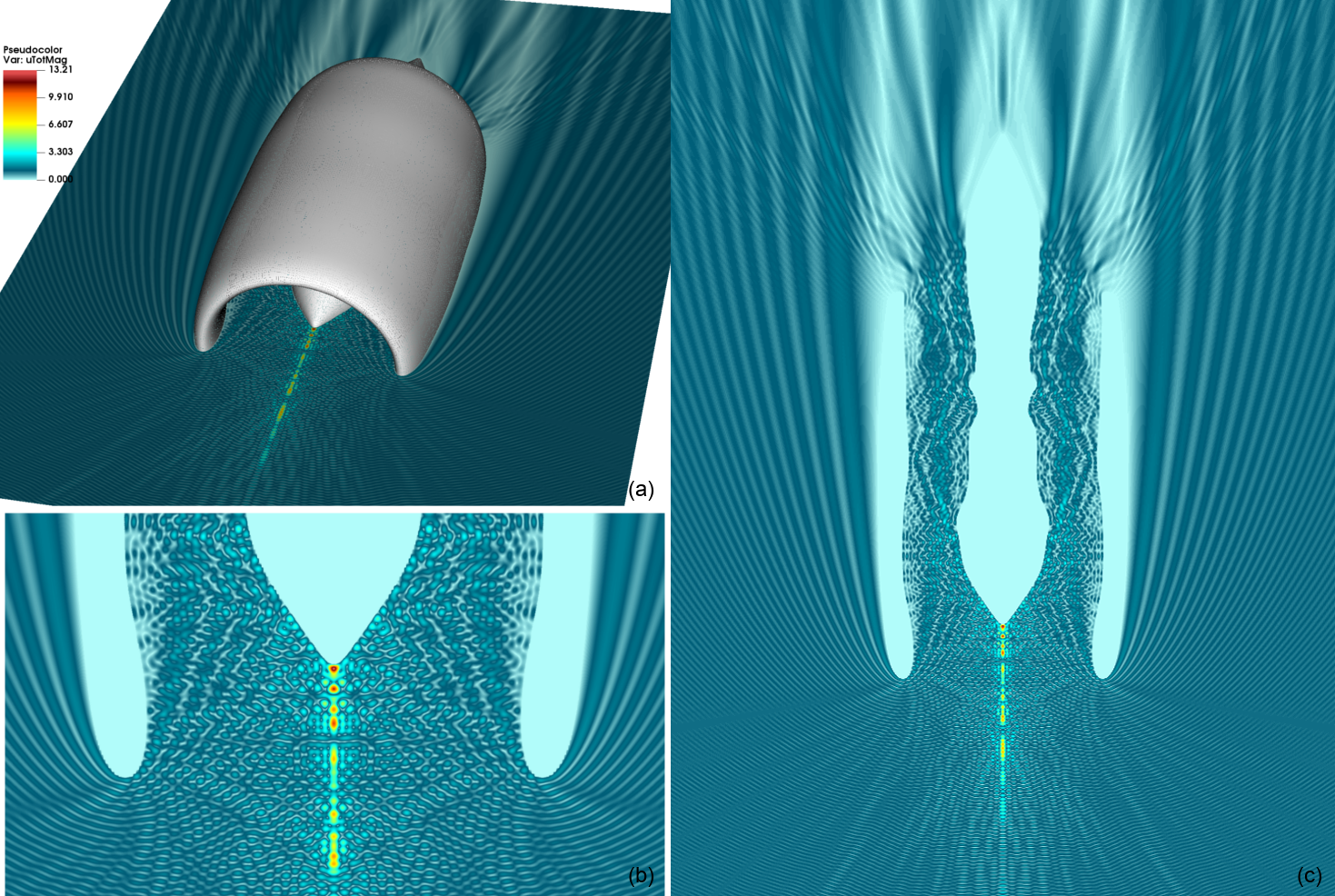}
    \caption{\small Total field magnitude $|u(x)| = |u^i(x) + u^s(x)|$
      pseudo-color plots for an $82$-wavelength aircraft nacelle. The
      field is plotted over a uniform grid of $2800\times 4000$ points
      $(x,z) \in [-4,4]\times[-5,6]$.  The nacelle is aligned with the
      $z$-axis and its front points towards the positive
      $z$-direction. The incident plane wave impinges on the geometry
      head-on, traveling towards the negative $z$-axis.}
  \label{fig:NacNFPlWav}
\end{figure}
\begin{table}[htb!] 
\centering
\begin{tabular}{lccccccc}\toprule
  & & & & \multicolumn{2}{c}{Plane Wave Incidence} & \multicolumn{2}{c}{Point Source Incidence} 
  \\ \cmidrule(lr){5-6} \cmidrule(lr){7-8}
  Unknowns   &Size         &IFGF levels  &Time (1 iter.)  
                     &Tot. iter. &$\nearErr$ &Tot. iter. &$\nearErr$  \\ 
  \midrule
  77,184     &10.2$\lambda$  &6            &2.0 s           
                     &33         & $1.5\cdot 10^{-3}$   &39        & $4.1\cdot 10^{-3}$ \\
  308,736    &20.5$\lambda$  &7            &8.7 s          
                     &47         & $2.2\cdot 10^{-3}$   &59        & $4.0\cdot 10^{-3}$ \\
  1,234,944  &40.9$\lambda$  &8            &40.4 s         
                     &55         & $6.7\cdot 10^{-4}$  &115       & $2.4\cdot 10^{-3}$ \\
  4,939,776  &81.8$\lambda$  &9            &176.4 s         
                     &65         & $1.5\cdot 10^{-3}$ (est.)   &219       & $3.5\cdot 10^{-3}$ (est.) \\
 \bottomrule
\end{tabular}
\caption{\small Near field relative errors $\nearErr$ resulting from
  applications of the IFGF accelerated acoustic solver to the nacelle
  geometry depicted in Figure~\ref{fig:NacelleAndMesh}, for plane wave
  and multiple point-source incidence, and with acoustical sizes
  ranging from $10.2$ to $81.9$ wavelengths. As detailed in the text,
  errors were evaluated by means of convergence studies, except for
  the $81.9$-wavelength case, for which an estimate was used.  The
  table additionally summarizes the total number of surface unknowns,
  the nacelle size in wavelengths, the maximum number of IFGF levels,
  the time required per GMRES iteration, and the total number of
  iterations incurred. In all cases the residual tolerance for the
  integral density was set to $10^{-3}$.}
\label{table:Nacelle}
\end{table}

Table~\ref{table:Nacelle} presents near field errors $\nearErr$
resulting from use of the proposed accelerated solver for an aircraft
engine nacelles $10.2,\, 20.5$, $40.9$ and $81.8$ wavelengths in size
under both plane wave and multiple point-source incidence.  The number
of IFGF levels is selected so that the finest-level IFGF box side
length is approximately $0.6\lambda$ in all cases, and the GMRES
residual tolerance was set to $10^{-3}$. The total near field relative
error $\nearErr$ was evaluated over a near field planar grid
$\Pnear^{xy}(z_0)$, where
$[x_{min},x_{max}] \times [y_{min},y_{max}] = [-4,4]^2, \,z_0 = -5$,
and $N_x = N_y = 400$. The reference solution used for error
estimation was obtained with the same number of surface patches but
using $8 \times 8$ points per patch and a residual tolerance of
$10^{-5}$, which gave rise to an increase in the number of iterations
required for convergence to the prescribed tolerance by only $8\--14$
iterations. In addition to the near field relative error,
Table~\ref{table:Nacelle} also includes the total number of unknowns,
the nacelle size in wavelengths, the time required per iteration, and
the total number of GMRES iterations required to meet the $10^{-3}$
residual tolerance. Note that, as the problem size increases from
$10.2\lambda$ to $20.5\lambda$, $20.5\lambda$ to $40.9\lambda$, and
$40.9\lambda$ to $81.8\lambda$, and the number of unknowns is
quadrupled in each case, the computing cost per iteration increases by
a factor of only $4.4,\,4.6$ and $4.4$, respectively (which is
consistent with an $\Ord(N\log N)$ complexity), and not the 16-fold
cost increase per wavelength doubling that would result from a
non-accelerated algorithm with quadratic complexity. This scaling of
the IFGF-accelerated combined-layer solver is also consistent with the
corresponding scalings presented in~\cite{BauingerBruno2021}, which do
not include singular local interactions, and suggests that the
partitioning and discretization of the geometry works well in
conjunction with the IFGF algorithm.  The results also indicate that
the discretization and $10^{-3}$ residual tolerance is sufficient to
produce solutions for the $10.2,\,20.5$ and $40.9$ wavelength cases
with an average error of $1.5\cdot 10^{-3}$ for plane wave scattering
and $3.5\cdot 10^{-3}$ for point source scattering. The average
relative error values were used to estimate the accuracy of the
$81.8\lambda$ simulation, which, like the smaller problems, converged
to the prescribed GMRES tolerance.

The total near field magnitude $|u(x)| = |u^i(x) + u^s(x)|$ for the
$81.8$-wavelength plane wave scattering case is displayed in
Figure~\ref{fig:NacNFPlWav}. The field magnitude is presented over the
$xz$-planar grid $\Pnear^{xz}(y^0)$ (equation~\eqref{eq:PlanarGrid}),
with
$[x_{min},x_{max}] \times [z_{min},z_{max}] = [-4,4] \times [-5,6],
\,y^0 = 0$, and $N_x = 2800$ and $N_z = 4000$. Along most of the
exterior circumference of the nacelle housing, the total field forms a
relatively uniform stratified pattern.  In other regions, intricate
multiple-scattering patterns develop, particularly in the region
around the intake and throughout the inside of the nacelle.
Figures~\ref{fig:NacNFPlWav}(b) and~(c) present views of the field
from the top, but with the scattering surfaces removed. Clearly, the
strongest field concentrations occur directly in front of the tip of
the nacelle shaft.  Note the symmetry in the detail of the near field
presented in Figure~\ref{fig:NacNFPlWav}(b), which results in spite of
the lack of symmetry, noted above in this section, in the geometry
discretization illustrated in the inset in
Figure~\ref{fig:NacelleAndMesh}(c)---which provides an additional
indication of the accuracy of the solution displayed.

\begin{figure}[h!]
    \centering
    \includegraphics[width=\textwidth]{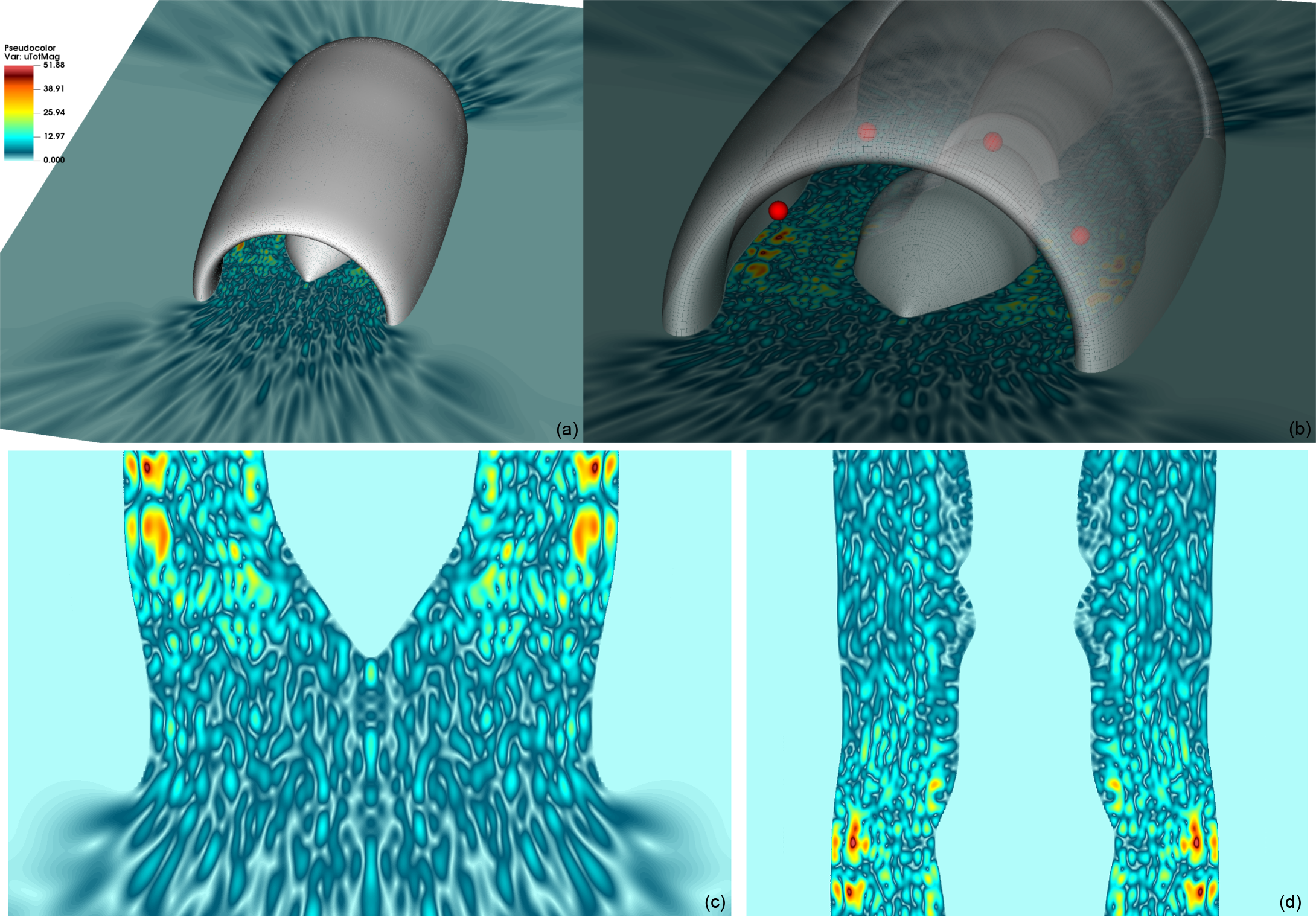}
    \caption{\small Total field magnitude $|u(x)| = |u^i(x) + u^s(x)|$
      pseudo-color plots for an $82$-wavelength aircraft nacelle. The
      field is presented over a uniform grid of $2800\times 4000$
      points $(x,z) \in [-4,4]\times[-5,6]$.  The nacelle is aligned
      with the $z$-axis and its front points towards the positive
      $z$-direction. The incident field is given by the
      sum~\eqref{eq:PointSource} of eight point sources within the
      nacelle around the center shaft, four of which are displayed as
      small red spheres in panel~(b).}
  \label{fig:NacNFPtSrc}
\end{figure}
Near fields resulting for the eight-point source/$81.8\lambda$ problem
are presented in Figure~\ref{fig:NacNFPtSrc}.  The total field
magnitude is presented over the $xz$-planar grid $\Pnear^{xz}(y^0)$
used previously in this section.  Figure~\ref{fig:NacNFPtSrc}(a)
presents the fields resulting under point-source incidence; note from
this image that the resulting fields scatter and exit the front inlet
and rear exhaust.  The close-up view in Figure~\ref{fig:NacNFPtSrc}(b)
highlights the location of four of the eight point sources, drawn as
red spheres for emphasis; the remaining four sources are obstructed
from view by the near field plane.  In Figures(c) and~(d) the geometry
is not included, enabling examination the field interaction within the
scatterer in greater detail. These two images exhibit complex multiple
scattering and, as mentioned in the context of the plane-wave
incidence problem considered earlier in this section, a high degree of
symmetry throughout the interior of the structure and in the regions
outside that surround the nacelle assembly. In contrast to the plane
wave scattering case, where the incident wave travels in directions
mostly parallel to the housing and shaft, placing sources between the
shaft and nacelle walls guarantees that most waves are scattered many
times by the surface before exiting the geometry.

\begin{figure}[h!]
    \centering
    \includegraphics[width=\textwidth]{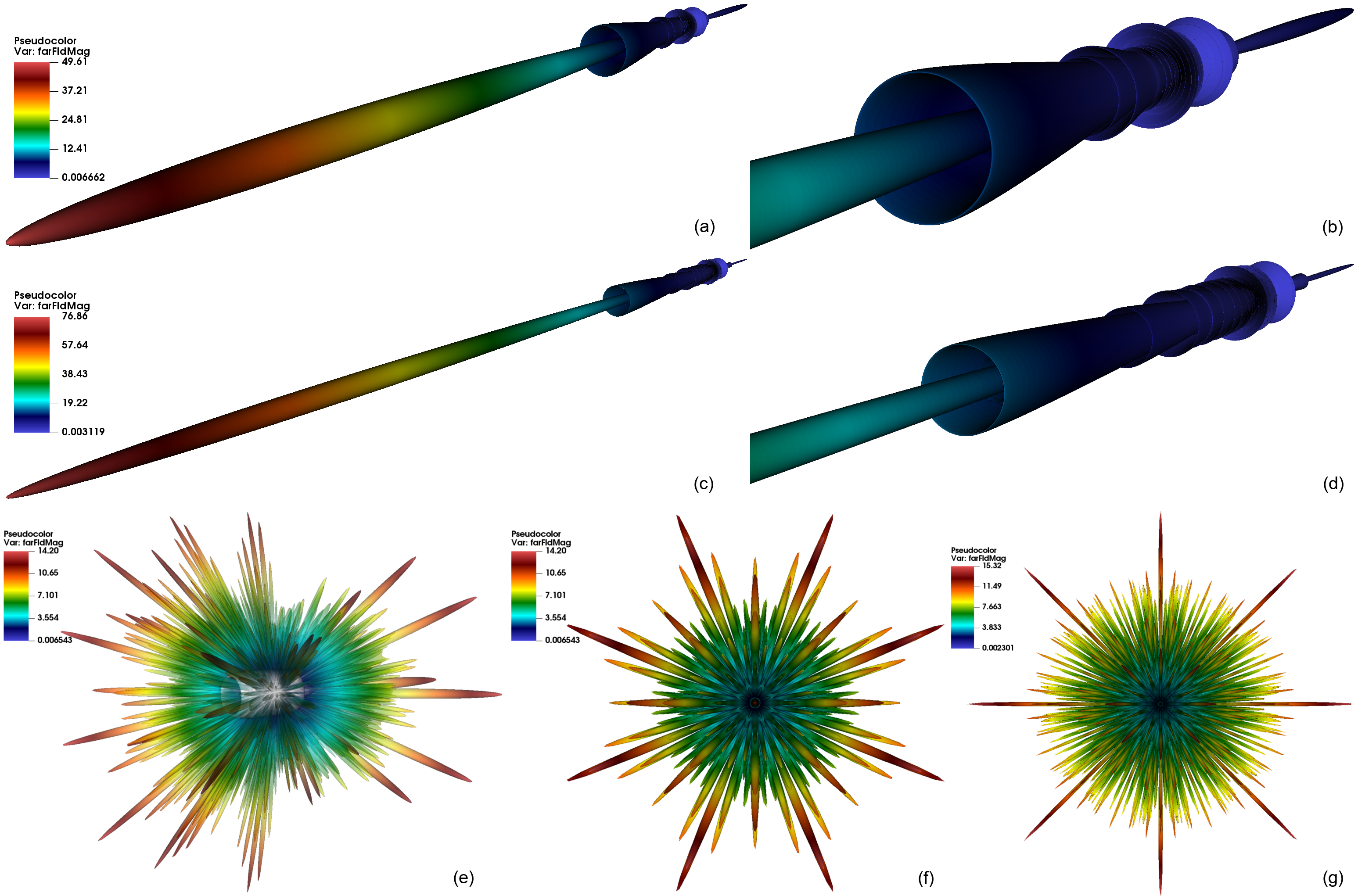}
    \caption{\small Far-field magnitudes for the nacelle geometry under plane
      wave and eight point-source incident fields. Panels (a) and~(b)
      present the far field for a $40.9\lambda$ plane wave and panels (c)
      and~(d) present the far field for an $81.8$-wavelength plane
      wave. Panels~(e)-(g) display far fields for the eight-point
      source incident field defined in~\eqref{eq:PointSource}, for
      $40.9\lambda$ in panels~(e) and~(f) and for $81.8\lambda$ in
      panel~(g).}
  \label{fig:NacelleFF}
\end{figure}
Figure~\ref{fig:NacelleFF} displays far-field magnitudes for both the
plane-wave and multiple-point source problems considered in this
section. In detail, Figures~\ref{fig:NacelleFF}(a)
and~\ref{fig:NacelleFF}(b) present the far field for a $40.9\lambda$
plane wave; and Figures~\ref{fig:NacelleFF}(c)
and~\ref{fig:NacelleFF}(d) present the far field for an
$81.8$-wavelength plane wave, all of which were computed using the
relation~\eqref{eq:FarField} over the discrete spherical
discretization $\Sfar$ (equation~\eqref{eq:SphereGrid}) with
$(N_{\phi},N_{\theta}) = (2000,600)$ in the $40.9\lambda$ case, and
with $(N_{\phi},N_{\theta}) = (3200,800)$ in the $81.8\lambda$ case.
These far field plots once again highlight the fact that most of the
wave energy is reflected toward the region directly in front of the
nacelle intake; note, also, that this reflection increases in
intensity as the frequency increases: the maximum far-field magnitude
increases by a factor of approximately $1.5$ as the problem size is
increased from $40.9\lambda$ to $81.8\lambda$.
Figures~\ref{fig:NacelleFF}(e-g) show the far field magnitude for the
eight-point source incident field~\eqref{eq:PointSource} at $40.9$ and
$81.8$ wavelengths. Figure~\ref{fig:NacelleFF}(e) displays a side view
of the $40.9\lambda$ far field magnitude $|\tilde{u}^{\infty}|$
including the nacelle geometry, for reference, with the intake
pointing left. Figures~\ref{fig:NacelleFF}(f)
and~\ref{fig:NacelleFF}(g), where the geometry is not included,
present the far field, with the positive $z$ direction pointing out of
the page, for the $40.9\lambda$ and $81.8\lambda$ cases, respectively.

\section{Conclusions}\label{sec:Conclusion}

The IFGF/RP scattering solver introduced in this paper demonstrates,
for the first time, the character of the IFGF acceleration algorithm
in the context of an actual problem of scattering. The combined
approach achieves accurate scattering simulations on the basis of an
$\Ord(N\log N)$ computational cost. Without recourse to the FFT
algorithm, and relying, instead, on a recursive interpolation scheme,
the IFGF accelerator enables the fast evaluation of the slowly-varying
factors in the factored Green function representations for groups of
sources. The IFGF algorithm is additionally integrated seamlessly with
the high-order RP methodology.  Two different strategies are proposed
in Section~\ref{sec:IFGF} for the treatment of the double layer
operator, each one of which, for efficiency, is recommended for use
under certain portions of the IFGF recursion. The OpenMP
implementation reported in this paper delivers accurate results in
brief computing times, and with fully converged solutions within the
accuracy quoted, even for the relatively large acoustical sizes and
for the complex and realistic scattering structures
considered. Comparisons with recent FMM-based solvers, presented in
Section~\ref{compari}, demonstrate the benefits provided by the
proposed IFGF/RP approach.

Beyond providing simulations of acoustic scattering by a sphere of up
to $128$ wavelengths in diameter, this paper demonstrated the
versatility of the proposed algorithms by including computational
results for two realistic engineering geometries, a submarine and a
turbofan nacelle.  The present contribution only considered a
shared-memory implementation of acoustic IFGF-based integral equation
solvers; the numerical examples presented show that even in this case
the proposed algorithm enables the efficient solution of problems with
millions of unknowns over complex geometries on the basis of a small
computer system. Owing in part to its virtually unique ability to
efficiently accelerate scattering operators without recourse to the
FFT, the IFGF algorithm can effectively be parallelized under the MPI
interface as well~\cite{BauingerBruno_parallel_2021}---and, thus, we
believe the IFGF/RP algorithm could advantageously be implemented in
distributed-memory systems. The development of such a MPI/OpenMP
parallel IFGF/RP implementation, which lies beyond the scope of this
paper, is left for future work.

\section*{CRediT authorship contribution statement}
\textbf{Oscar P. Bruno:} Conceptualization, Methodology, Validation,
Investigation, Resources, Writing, Supervision, Funding acquisition.
\textbf{Edwin Jimenez:} Conceptualization, Methodology, Software, Validation,
Investigation, Writing, Visualization.  \textbf{Christoph Bauinger:}
Conceptualization, Methodology, Software, Validation, Investigation, Writing,
Visualization.

\section*{Declaration of competing interest}
The authors declare that they have no known competing financial interests or
personal relationships that could have appeared to influence the work reported
in this paper.

\section*{Acknowledgments}

This work was supported by NSF, DARPA and AFOSR through contracts
DMS-2109831 and HR00111720035 and FA9550-21-1-0373, and by the NSSEFF
Vannevar Bush Fellowship under contract number N00014-16-1-2808.

\bibliographystyle{unsrt}
\bibliography{main}



\end{document}